\documentclass[a4paper, 
10pt
]{article}

\usepackage[a4paper,left=3cm,right=3cm,top=2.5cm,bottom=2.5cm]{geometry}
\usepackage{amsmath}
\usepackage{multirow}
\usepackage[table]{xcolor}

\usepackage{subcaption}
\usepackage{pgfplots}
\usepackage{pgfplotstable}
\usepgfplotslibrary{groupplots}
\usepgfplotslibrary{fillbetween}
\usepackage{mathtools}
\usepackage{multicol}
\usepackage{comment}
\usepackage{booktabs}
\pgfplotsset{compat=1.5}
\usepackage{amssymb}
\usepackage{url}
\usepackage{bm}

\usepackage{float}
\usepackage{tabularx}
\usepackage{enumerate}
\usepackage{array}
\usepackage{xspace}
\usepackage{tikz}
\usepackage{tikzsymbols}
\usetikzlibrary{calc,trees,positioning,arrows,chains,shapes.geometric,%
    decorations.pathreplacing,decorations.pathmorphing,shapes,%
    matrix,shapes.symbols, decorations.markings, patterns,fit,quotes,%
    arrows.meta}

\usepackage{authblk}
\usepackage{xcolor}
\definecolor{Gray}{gray}{0.9}

\usepackage{hyperref}
\usepackage{siunitx}

\usepackage[bottom]{footmisc}

\newcommand{\mupar}{\ensuremath{\boldsymbol{\mu}}}
\newcommand{\etapar}{\ensuremath{\boldsymbol{\eta}}}
\newcommand{\R}{\ensuremath{\mathbb{R}}}
\newcommand{\fulldim}{\ensuremath{\mathcal{N}}}
\newcommand{\reddim}{\ensuremath{N}}

\newcommand{\fullspace}{\ensuremath{\mathbb{X}^\fulldim}}

\newcommand{\parspace}{\ensuremath{{\bf P}}}

\newcommand{\RA}[1]{{\color{black}#1}}
\newcommand{\RB}[1]{{\color{black}#1}}
\newcommand{\RC}[1]{{\color{black}#1}}

\begin{document}

\title{Hull shape design optimization with parameter space and model reductions, and self-learning mesh morphing} 

\author[]{Nicola~Demo\footnote{nicola.demo@sissa.it}}
\author[]{Marco~Tezzele\footnote{marco.tezzele@sissa.it}}
\author[]{Andrea~Mola\footnote{andrea.mola@sissa.it}}
\author[]{Gianluigi~Rozza\footnote{gianluigi.rozza@sissa.it}}

\affil{Mathematics Area, mathLab, SISSA, via Bonomea 265, I-34136
  Trieste, Italy}

\maketitle

\begin{abstract}
In the field of parametric partial differential equations, shape
optimization represents a challenging problem due to the required computational
resources. In this contribution, a data-driven framework involving multiple
reduction techniques is proposed to reduce such computational burden. Proper
orthogonal decomposition (POD) and active subspace genetic algorithm (ASGA) are
applied for a dimensional reduction of the original (high fidelity) model and
for an efficient genetic optimization based on active subspace property.  The
parameterization of the shape is applied directly to the computational mesh,
propagating the generic deformation map applied to the surface (of the object to
optimize) to the mesh nodes using a radial basis function (RBF) interpolation.
Thus, topology and quality of the original mesh are preserved, enabling
application of POD-based reduced order modeling techniques, and avoiding the
necessity of additional meshing steps. Model order reduction is
  performed coupling POD and Gaussian process regression (GPR) in a
  data-driven fashion. The framework is validated on a
benchmark ship.
\end{abstract}

\tableofcontents

\section{Introduction}
\label{sec:intro}
In the framework of parameterized partial differential equation (PDE) problems
for engineering, reduced order models (ROMs) and optimization algorithms are two
instruments that particularly benefit a synergic use. In several cases
of engineering interest in which PDEs solution require considerable 
computational effort, ROMs enable in fact a remarkable reduction in the resources
required for each calculation. There are of course several ways to
reduce the dimensionality of discretized PDEs. The most naive approaches,
such as coarsening the computational grids clearly have negative effects
on the quality of the solutions. This is particularly true for
problems characterized by complex physics and geometrical features, which
in most cases require a very high number of degrees of freedom, ultimately
resulting in expensive computations. In the context of an optimization algorithm
execution, where many discretized PDE solutions must be computed, the
overall computational load often becomes unaffordable.
With only modest negative effects on the PDE solution accuracy, ROMs
can be conveniently exploited to reduce the high dimensionality of the
original discrete problem --- to which we will herein refer to as
full order model (FOM) or \textit{high fidelity} model.
ROM algorithms can be employed in several industrial design processes, and in
particular to shape optimization, in which the objective of the computations
is \RC{to find the best shape} of a particular product or artifact. Such problems are
in fact typically modeled through parametric PDEs, in which input parameters
control the geometric features of the object at hand. ROMs efficiently approximate
the numerical solution of the full order PDE with a suitable reduced surrogate,
enabling drastic reduction in the computational burden of the overall optimization
procedure.

There are of course several different algorithms which allow for an efficient
reduction of the dimensionality of parametric problem. In the present contribution,
we make use of a data-driven approach based on proper orthogonal
decomposition (POD)~\cite{salmoiraghi2016advances,
rozza2018advances}.
The equation-free nature of such method is often an essential feature in the
industrial sector, where modularity and solvers \textit{encapsulation} play a
fundamental role. Indeed, the data-driven POD based ROM employed in the present
optimization framework can be coupled with any PDE solver, as the data integration
is enforced through the output of interest of the full order problem.
Similar reduced methods have been proposed
in~\cite{demo2018shape,DemoOrtaliGustinRozzaLavini2020BUMI} for the shape
optimization of a benchmark hull, while additional improvements have been made
coupling the ROM with active subspace analysis and different shape parameterization algorithms
in~\cite{tezzele2018dimension,demo2018isope,demo2019marine,tezzele2019marine}. We refer
the readers interested in parametric hull shape variations using ROMs
to~\cite{villa2020parametric}, while we mention~\cite{diez2015design, serani2017towards}
for design-space dimensionality reduction in shape optimization with POD. 
Moving from hulls to propellers, data-driven POD has also been successfully
incorporated in the study of marine propellers efficiency~\cite{mola2019marine, gaggero2020reduced}
as well as hydroacoustics performance~\cite{gadalla2020les}.

A further aspect of novelty of the optimization framework proposed is
\RC{related to}
the parameterization of the geometry. In typical shape optimization cycles, the
surface of the object under study is deformed before the domain discretization
takes place. Thus, the meshing phase is repeated for any deformed entity. Such
approach has the clear advantage of allowing for good control of the
quality of the computational grid produced for each geometry tested. Yet, it
suffers of two main problems: \textit{i)} the meshing
step may be expensive, both because its CPU time might be comparable to the
resolution of the problem itself, and because mesh generation is specially intensive
in terms of human operator hours required; \textit{ii)} a different mesh 
for each geometry does not allow for the application of POD or several other
ROM approaches, which require that the mesh topology, as well as the number
of degrees of freedom of the discretized problem, are conserved across all
the shapes tested. Thus, assuming a generic deformation map is available,
which morphs the initial object surface --- not the grid ---, we exploit such
deformation to train a radial basis function (RBF) interpolation that will extend
the surface deformation to the nodes of the PDE volumetric mesh. In
this sense, the method is capable to learn and propagate any
deformation to a given mesh. Properly selecting 
the RBF kernel, we can then obtain a smooth deformation in all the discretized domain,
not only ensuring that the overall parameterization map preserves the initial mesh quality
but also its topology. We remark that in this work, free-form deformation (FFD) is used
to deform the surface of the object under study. Yet, we stress that the RBF extension
methodology is completely independent from the parameterization method chosen for the
object geometry. \RB{A similar approach has been recently investigated in~\cite{jeong2020mesh}.}

The optimization algorithm used in this work is the recently developed
active subspaces extension of the classical genetic algorithm called
ASGA~\cite{demo2020asga}, which performs the mutation and cross-over steps on a
reduced dimensional space for a faster convergence.

\begin{figure}[tb]
\centering
\includegraphics[width=0.8\textwidth]{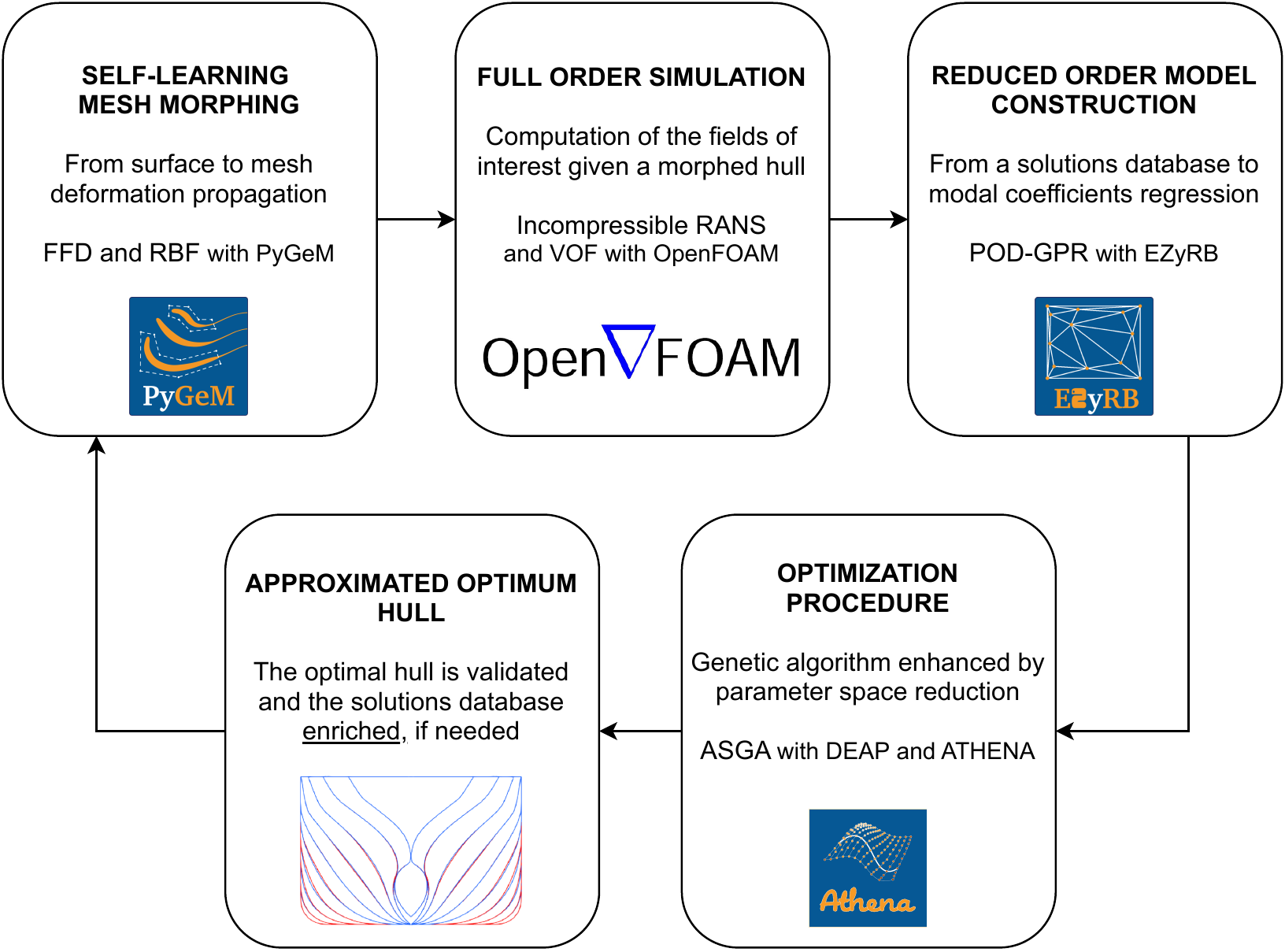}
\caption{\RB{Illustration of the key steps of the proposed optimization
pipeline with the methods and the softwares used.}}
\label{fig:paper_scheme}
\end{figure}

All the algorithms used in this work are implemented in open source software
libraries~\cite{TezzeleDemoMolaRozza2020PyGeM,romor2020athena,DemoTezzeleRozza2018EZyRB,gpy},
which we will briefly introduce in the discussions of the
corresponding numerical methods. \RB{In Figure~\ref{fig:paper_scheme}
we depicted an outline of the whole numerical pipeline we are going to
present, emphasizing the methods and the softwares used.}
One of the main goals of this contribution it that of testing the full pipeline composed by
data-driven POD ROM, combined FFD-RBF shape parameterization algorithm and ASGA optimizer
on a problem that can be both meaningful to the ship hydrodynamics community and
easily reproducible. For such reason, the test case considered is that of the DTC
hull~\cite{moctar2012duisburg}, for which online tutorials are available to run fairly
accurate flow simulations in fixed sink and trim conditions. Since in such set up, the
hull optimizing resistance is a trivial, zero volume hull, the DTC benchmark hull is
here optimized based on the \RB{total resistance coefficient  $C_t$}.
We organize the
contribution as follows: Section~\ref{sec:shape} presents a deeper discussion about
the parameterization of the object and of the computational grid;
Section~\ref{sec:model} describes the full order model and the reduced order
one, while Section~\ref{sec:asga} is devoted to an algorithmic discussion about
the optimization algorithm and its supporting mathematical tools. The final
sections, \ref{sec:results} and \ref{sec:conclusions}, show the numerical
results obtained and present the conclusive summary, respectively.

\section{Shape and grid parameterization}
\label{sec:shape}
Whenever industrial design processes as the ones discussed in this work
are aimed at improving, among other aspects, the geometric features of a particular
artifact, a shape parameterization algorithm is a cornerstone of the whole optimization
pipeline. Optimization tools, as well as the non-intrusive model
reduction techniques employed in the present investigation, are in fact based on the parameterized PDEs
paradigm introduced in the previous section. In such framework, a set of geometric
input parameters affects the output of a parametric PDE through the deformation of its
domain geometry. Thus, the shape parameterization algorithm role is that of
mapping the variation of a set of numerical parameters, to the corresponding
deformation of the PDE domain geometry. In other words, since optimization
tools are mathematical algorithms which must be fed with numbers, the
shape parameterization algorithms translate shape deformations into
variations of the numeric quantities they need.

\subsection{How to combine different shape parametrization strategies}
\begin{figure}[t]
\centering
\includegraphics[height=0.14\textheight]{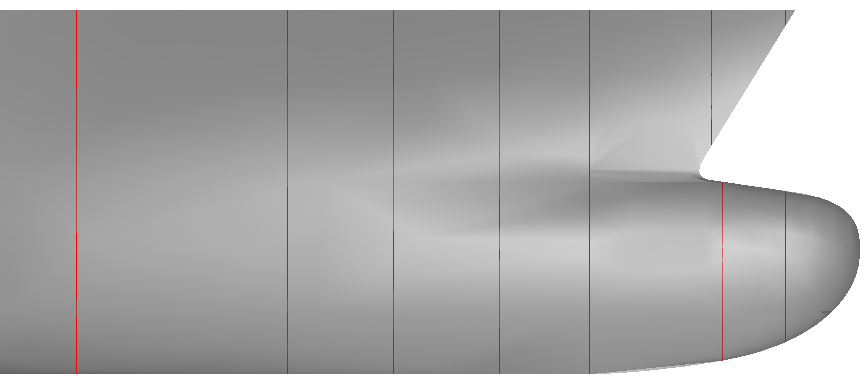}
\includegraphics[height=0.14\textheight]{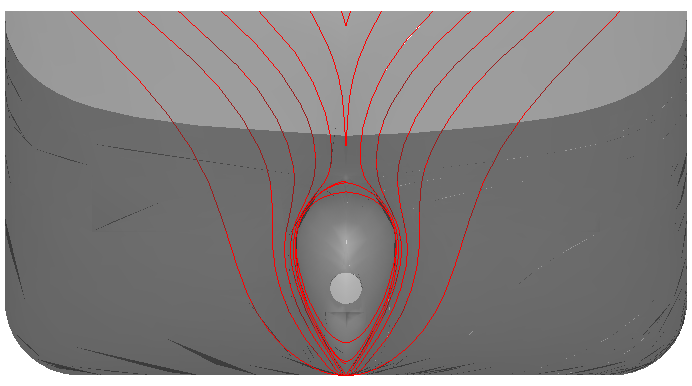}
\caption{Side view (left) and front view (right) of a typical container ship hull
         bow region.}
\label{fig:KCS_bow}
\end{figure}

In this work, we make combined use of two general purpose shape
parameterization algorithms to deform the three dimensional geometry of a ship hull,
and accordingly update the volumetric grid used for ship hydrodynamics simulations
in a fully automated fashion. More specifically, free form deformation (FFD) is first
used to generate a family of deformations of the surface of a base hull. In a second
step, radial basis functions (RBF) interpolation is used to propagate the hull surface deformation to
the internal nodes of the fluid dynamic simulation computational grid. For visual reference,
Figure~\ref{fig:KCS_bow} depicts the side view (on the left) and front view (on the right)
of a container ship hull bow region. In the picture, several sections perpendicular to
the hull longitudinal axis are indicated by red lines.

Despite an extensive discussion of FFD and RBF theoretical foundations is clearly
beyond the scope of the present contribution, this section will introduce the
key concept upon which both algorithms are based and describe their combined
deployment in the framework of our optimization pipeline.

\begin{figure}[htb]
\centering
\resizebox{\textwidth}{!}{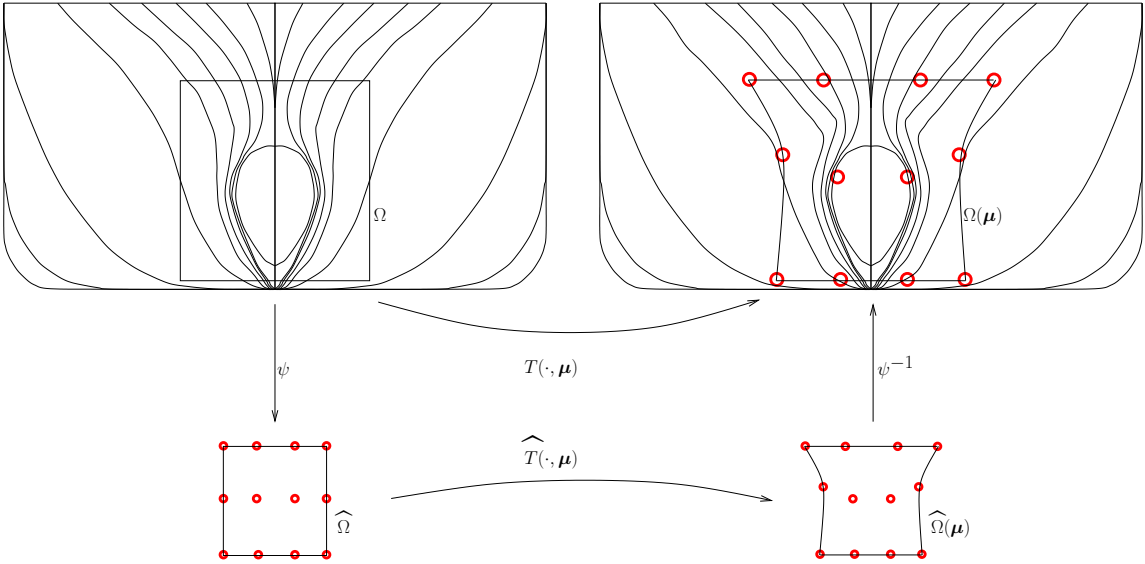}
\caption{A two dimensional sketch of the FFD procedure applied to the surface of
         a container ship hull\RB{, including the three transformations
         $\psi$, $\widehat{T}(\cdot,\mupar)$ and $\psi^{-1}$
         composing the process}.}
\label{fig:FFD_sketch}
\end{figure}

The first shape parameterization algorithm applied in this work is
the free form deformation~\cite{sederbergparry1986,LassilaRozza2010,sieger2015shape}.
As mentioned, it is a general purpose algorithm, designed to be applied
to arbitrarily shaped geometries. FFD is fundamentally made up of three different
geometrical transformations, as illustrated in Figure~\ref{fig:FFD_sketch}.
The first transformation $\psi$ maps the physical domain $\Omega$ into a reference
domain $\widehat{\Omega}$.
In such domain, a lattice of points is generated, and are used as the control points
of a set of smooth shape functions such as the Bernstein polynomials used in this work.
Thus, once a displacement is prescribed to one or more of the control points in the
lattice, the shape functions are used to propagate such displacement to all the points
in the reference domain $\Omega$. The smooth displacement field obtained, is
the second and most important transformation $\widehat{T}$ in the FFD process.
In the third, final step, the deformed reference domain is mapped back into
the physical one by means of $\psi^{-1}$ to obtain the resulting morphed geometry.

The current description suggests that the parameters $\mupar$
of the final FFD map $T(\cdot,\mupar)$ are the displacements prescribed
to one or more of the lattice control points. The procedure can account for
both a variable number of lattice points and of displaced control points.
For such reason, FFD deformations can be built with an arbitrary number
of parameters. 

We point out that the FFD algorithm results in a displacement law
for each 3D space point within the control points lattice. Thus, it
can be readily deployed to deform shapes specified through
surface triangulations (such as STL geometries) and surface grids
in general. In addition, it can be also used to directly deform
volumetric grids used for fluid dynamic simulations. Yet, mainly for
practical reasons, in this work we only make
use of FFD to deform the STL surface triangulation describing the
hull geometry. In fact, we must point out that if FFD has to be used
to modify the volumetric mesh used for CFD simulations, the control
points lattice dimensions must be much bigger than those needed
when only deforming the hull surface, leading to infeasible
optimization procedures. This is due to the fact that
when deforming volumetric meshes, it is often convenient to distribute the
deformations over a high number of cells, rather than concentrating
all the displacements in a very confined region in which cells can get
distorted or even inverted. But because FFD only affects points located
within the control points lattice, this means that the latter must
extend for a bigger volume. In addition, to maximize the volumetric
mesh quality, the user must include more control points in the lattice
to make sure that different deformation magnitudes are imposed in
regions close to the hull and far from it. Such \emph{manual} control
over the local mesh deformation can often become quite cumbersome.

For such reasons, after the hull surface mesh has been modified by means
of FFD, we resort to RBF to propagate the hull boundary displacements to
the internal nodes of the volumetric mesh for CFD simulations. In a broader
sense, RBF is an interpolation algorithm, in which linear combinations of
radial bases are used to approximate a function with values prescribed
only in a finite number of points, in every point of a domain. In the case
of interest, the displacement field function prescribed on the points
of the hull surface must be interpolated in the positions corresponding to
every node of the volumetric mesh. Thus, the displacement obtained from the
$m$ surface nodes original position $\{\boldsymbol{s}_1,\dots,\boldsymbol{s}_m\}$
and the corresponding displaced position $\{\boldsymbol{s}'_1,\dots,\boldsymbol{s}'_m\}$
must be interpolated at the positions $\{\boldsymbol{v}_1,\dots,\boldsymbol{v}_n\}$
of the $n$ volumetric mesh nodes. Such interpolation reads
\begin{equation}
\label{eq:RBF_interp}
\boldsymbol{d}(\boldsymbol{x}) = \sum_{j=1}^m \boldsymbol{w}_j\varphi_j(\boldsymbol{x}),
\end{equation}
where the radial bases $\varphi_j(\boldsymbol{x}) = \varphi_j(||\boldsymbol{x}-\boldsymbol{x}_j||)$
are functions that only depend on the distance between evaluation point $\boldsymbol{x}$ and control point
$\boldsymbol{x}_j$.
The weights $\boldsymbol{w}_j$ are computed by imposing the interpolation constraints
$\boldsymbol{d}(\boldsymbol{s}_i)=\boldsymbol{s}'_i-\boldsymbol{s}_i$, after a
radial basis has been centered at every constrained point ($\boldsymbol{x}_j=\boldsymbol{s}_j$).
This results in the linear system
\begin{equation}
\label{eq:linSys}
\boldsymbol{A}\boldsymbol{X} = \boldsymbol{B},
\end{equation}
where
\begin{equation}
\label{eq:linSysMat}
\boldsymbol{A}=\left[
\begin{array}{c c c}
\varphi_1(\boldsymbol{s}_1) & \dots & \varphi_1(\boldsymbol{s}_m) \\
\vdots & \ddots & \vdots \\
\varphi_m(\boldsymbol{s}_1) & \dots & \varphi_m(\boldsymbol{s}_m) 
\end{array}
\right],\qquad
\boldsymbol{X} =\left\{
\begin{array}{c}
\boldsymbol{w}_1 \\
\vdots \\
\boldsymbol{w}_m 
\end{array}
\right\},\qquad
\boldsymbol{B} =\left\{
\begin{array}{c}
\boldsymbol{s}'_1-\boldsymbol{s}_1 \\
\vdots \\
\boldsymbol{s}'_m-\boldsymbol{s}_m
\end{array}
\right\}.
\end{equation}

Linear system \eqref{eq:linSys} is solved in a pre-processing phase, and the weights
computed are then used to compute the displacement of every node of the volumetric mesh
by means of Equation~\eqref{eq:RBF_interp}. The latter operation can be conveniently
carried out in a parallel fashion, and is highly efficient. On the other hand, 
$\boldsymbol{A}$ is a full $m\times m$ matrix which can make the solution of
system~\eqref{eq:linSys} quite time and memory demanding when a large number of
RBF control points are considered. That is why, in some cases only a portion of the
surface mesh nodes are used as RBF control points, which limits the computational cost
more than linearly, and in most cases has only modest effect on the morphing accuracy.
\begin{figure}[t]
\begin{tabular}{c c}
\includegraphics[height=0.2\textheight]{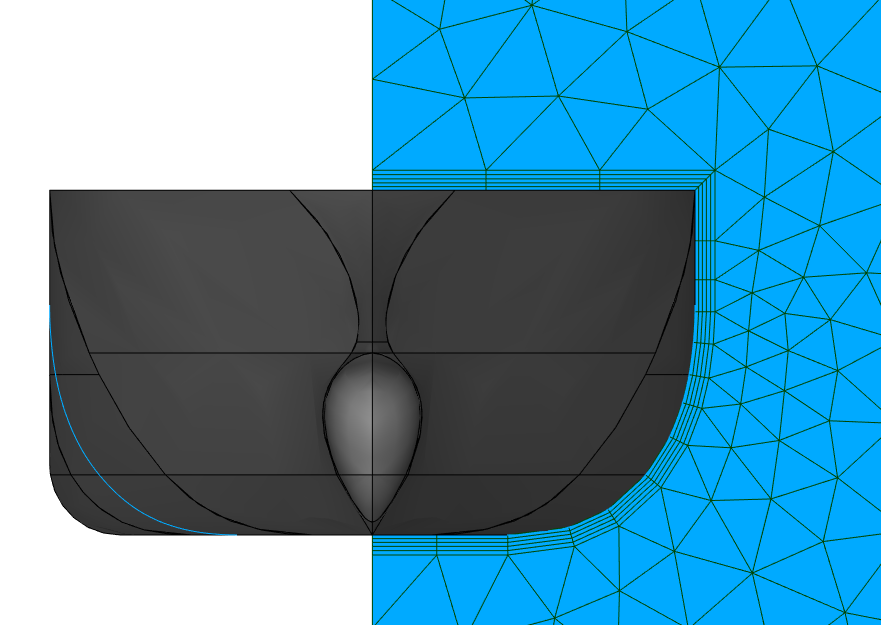} &
\includegraphics[height=0.2\textheight]{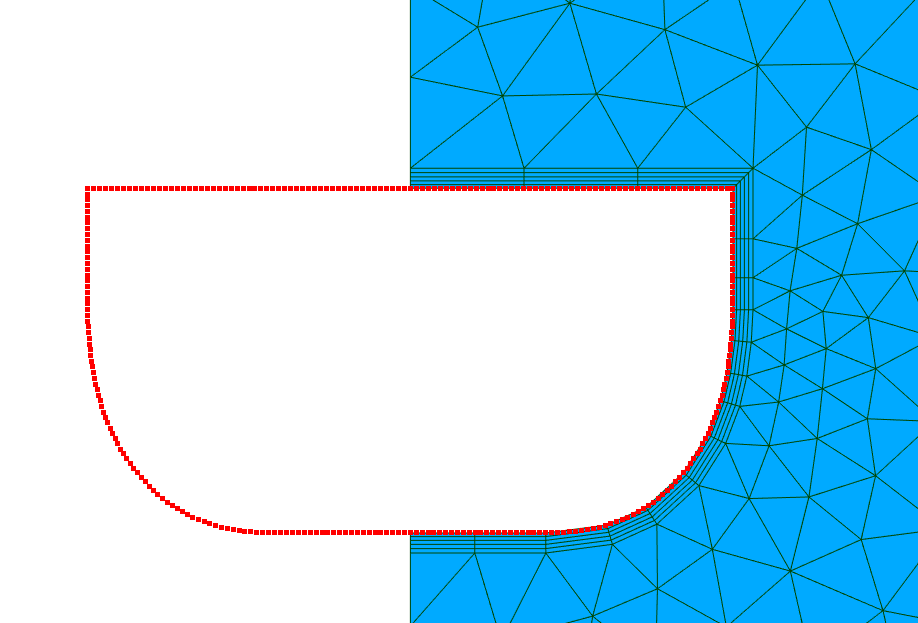} \\
a) & b) \\[0.6mm]
\hspace{-2mm}\includegraphics[height=0.2\textheight]{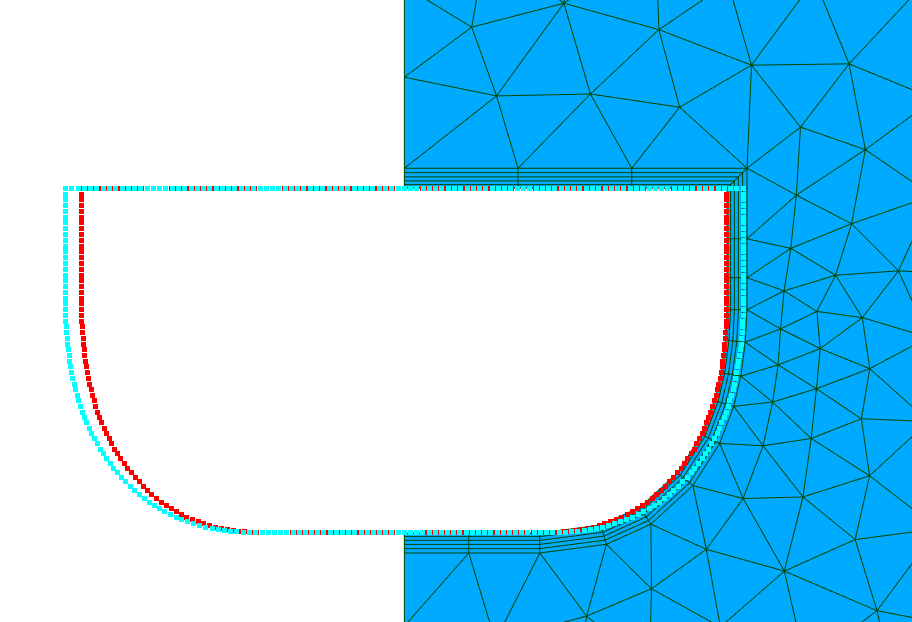} &
\includegraphics[height=0.2\textheight]{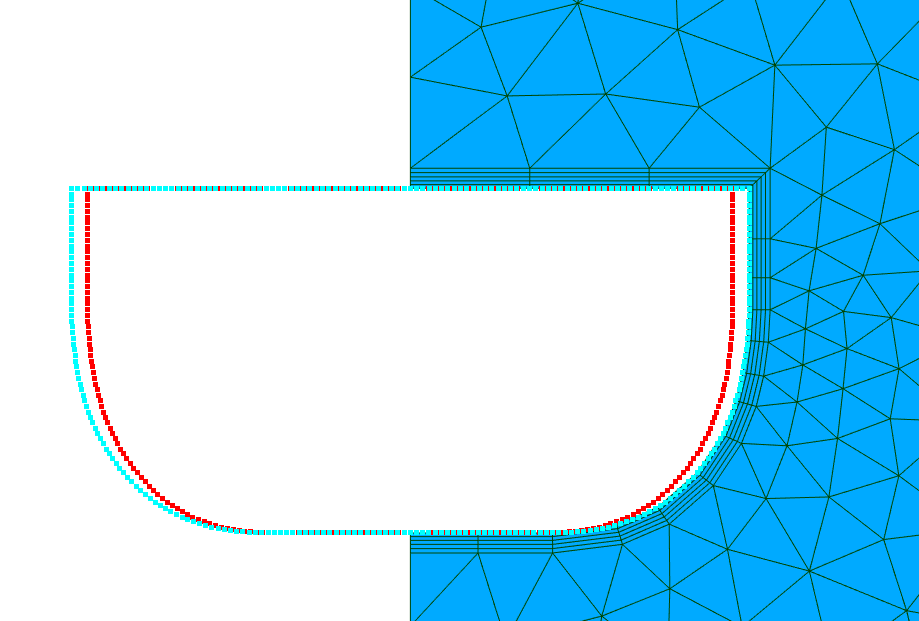} \\
c) & d) \\ \\
\end{tabular}
\caption{\RA{A section view example illustrating the RBF morphing steps carried out to
         propagate the hull surface deformations to a volumetric mesh for ship hydrodynamics
         simulations.}}
\label{fig:RBF_KCS}
\end{figure}

Both the FFD and RBF algorithms briefly described in this section have been implemented
in the Python library for geometrical morphing PyGeM~\cite{TezzeleDemoMolaRozza2020PyGeM},
which has been used to produce all the deformed geometries and computational grids used
in this work. An example of the RBF application to volumetric mesh morphing described in
this paragraph is presented in Figure~\ref{fig:RBF_KCS}. The figure illustrates all the
steps involved in the procedure, which starts with a) a first volumetric mesh around
the hull, and b) a surface mesh on the hull surface. In step c) the latter mesh is then
deformed and d) the surface mesh displacement field is finally used to feed the RBF algorithm
and propagate the boundary motion to the internal volumetric mesh nodes. As it can be appreciated
in the illustration, to avoid distortion of the volumetric mesh symmetry plane, the surface
mesh must include both sides of the hull. In the present work, the deformation of the surface
mesh has been carried out by means of FFD. Yet, we remark that any deformation law which 
results in a one to one correspondence between original and deformed surface grids can be
propagated to the nodes of the volumetric mesh with RBF interpolation.

\section{The mathematical model for incompressible fluids}
\label{sec:model}

The computational gain of the proposed pipeline is obtained by using a
model order reduction based on proper orthogonal decomposition (POD) to approximate the solution of the parametric
PDEs describing the studied phenomenon. This
technique assumes an initial solutions database produced by solving the full
order model (FOM), for some values of the parameters. We refer to such solutions as
high-fidelity solutions, or \textit{snapshots}. Depending on the intrusiveness
of the reduced order method, also the discrete operators of the numerical
problem can be required.
In this contribution, we propose a non-intrusive approach, constructing a reduced
order model (ROM) within a data driven setting using the FOM snapshots and the corresponding parameter
values (described in Section~\ref{sec:shape}). This allows a modular
structure where any numerical solver, also commercial, can be adopted,
since the ROM relies only on input and output couples. 

The following paragraphs present the full order model used in this
work and the ROM constructed with it. We briefly describe the incompressible Reynolds Averaged
Navier--Stokes (RANS) equations and its numerical solution in a finite volume
framework, then we proceed with an algorithmic analysis of the proper
orthogonal decomposition with Gaussian process regression (POD-GPR).

\subsection{The full order model: incompressible RANS}
The FOM used in this work is the Reynolds Averaged Navies--Stokes (RANS) model complemented by
a Volume of Fluid (VOF) front capturing method to deal with the multi phase nature of the fluid
surrounding the hull. The resulting govern equations are discretized by means of a Finite Volumes (FV)
strategy implemented in the open source library openFOAM~\cite{openfoam}. Such mathematical and
numerical setup is nowadays adopted in many industrial naval contexts thanks to its robustness and
accuracy. The test case considered is one of the tutorials of the library, which is designed to
reproduce the DTC experiments reported in reference~\cite{moctar2012duisburg}. We here
provide a minimal overall description of the model. We refer to the original documentation of the
library for all the numerical and technical details.

The RANS equations model the turbulent incompressible flow, while the volume of
fluid (VOF) technique~\cite{vof} is applied to handle the biphase nature of
the fluid (water and air).
The equations governing our system are the following
\RB{
\begin{equation}
\begin{cases}
\frac{\partial \bar{u}}{\partial t} + (\bar{u} \cdot \nabla)\bar{u}- \nabla \cdot
(\tilde{u} \otimes \tilde{u})
= -\frac{1}{\rho} \nabla \bar{p} +\nabla \cdot \nu \nabla \bar{u} +g,\\
\nabla \cdot \bar{u}=0,\\
\frac{\partial \alpha}{\partial t} + \nabla \cdot (\bar{u} \alpha)=0,\\
\end{cases}
\label{ns}
\end{equation}
}
where $\bar{u}$ and $\tilde{u}$ refer to the mean and fluctuating velocity after
the RANS decomposition, respectively, $\bar{p}$ denotes the mean pressure, $\rho$ is the
density, $\nu$ the kinematic viscosity, and $\alpha$ is the discontinuous
variable belonging to interval $[0, 1]$ representing the fraction of the second
flow in the infinitesimal volume. \RB{Finally, vector $g$ represents the body accelerations
associated with gravity.}

The first two equations are the continuity and momentum conservation, where the new term, the
Reynolds stresses tensor $\tilde{u} \otimes \tilde{u}$, have to be modeled with
additional equations in order to close the system.
Among all the turbulence models available in literature, we use the
$\text{SST} k-\omega$ turbulence model~\cite{menter}.
The third equation represents the transport of the VOF variable $\alpha$. 
Such variable controls also the density $\rho$ and the kinematic viscosity $\nu$, since they are
defined using an algebraic formula expressing them as a convex combination of
the corresponding properties of the two flows such that
\begin{equation}
\rho = \alpha \rho_{\text{air}} + (1-\alpha) \rho_{\text{water}},
\qquad \nu = \alpha \nu_{\text{air}} + (1-\alpha) \nu_{\text{water}}.
\end{equation}
To compute the steady solution in a discrete environment, we apply the finite volume
(FV) approach.  We set a pseudo--transient simulation, applying a first order
implicit local scheme for the temporal discretization, while for the
spatial scheme we apply the linear upwind one. Regarding the software,
as mentioned the
simulation is carried out using the C++ library OpenFOAM~\cite{openfoam}.

\subsection{The reduced order model: POD-GPR}
POD is a linear dimensional reduction technique capable to construct a reduced order
model from a set of high-fidelity snapshots.
Such space is spanned by (typically few) basis functions, that are computed by
minimizing the error between the original snapshots and their orthogonal
projection~\cite{volkwein-rom}.
In a parametric context, it enables --- provided a proper set of parameter
samples --- the possibility to approximate the solution manifold in a very
efficient way. Formally, we define the set of parameters
$\{\mupar_i\}_{i=1}^M$ such that $\mupar_i \in \parspace \subset \R^p$ for $i = 1,\dotsc, M$.
For each parameter, the solution is computed using the FOM. Let
$\fulldim$ be number of degrees of freedom of the full simulation, we obtain the solutions
$\mathbf{x}_i \in \fullspace_i$ for $i = 1, \dotsc, M$. Since the finite volume
space is created only once and then it is deformed, all the geometric
configurations have the same dimensionality even if they belong to different spaces.
The vectorial solutions are arranged as columns of the snapshots matrix, such that
\begin{equation}
\mathbf{X} = \begin{bmatrix} | & \dotsc & | \\ \mathbf{x}_1 & \dotsc & \mathbf{x}_M\\ | & \dotsc & | \\\end{bmatrix} \in \R^{\fulldim \times M}.
\end{equation}

The basis of the POD space, composed by the so called POD modes, is computed
using the singular value decomposition (SVD) of the snapshots matrix
$\mathbf{X} = \mathbf{U} \boldsymbol{\Sigma} \mathbf{V}^*$. 
The unitary matrix
$\mathbf{U} \in \R^{\fulldim \times M}$ contains the left-singular vectors of
$\mathbf{X}$, which are the POD modes. Moreover the diagonal matrix
$\boldsymbol{\Sigma} = \text{diag}(\lambda_1, \dotsc, \lambda_M)$, where
$\lambda_1 \ge \lambda_2 \ge \dotsc \ge \lambda_M$, contains the singular
values, which indicate the energetic contribution of the corresponding
modes. By looking at the spectral decay we can retain the first
$\reddim$ most energetic modes, which span the optimal space of dimension
$\reddim$.

Such basis can be exploited in a Galerkin projection
framework~\cite{StrazzulloBallarinRozza2020,girfoglio2020podgalerkin,rozza}
, in an
hybrid framework combining data-driven methods with
projection~\cite{HijaziStabileMolaRozza2020b,GeorgakaStabileStarRozzaBluck2020},
or used to project onto the reduced space the initial snapshots. Thus
we can approximate the snapshots $\mathbf{x}_j$ as a linear
combination of the modes as
\begin{equation}
\mathbf{x}_j = \sum_{i=1}^M \mathbf{c}^i_j \boldsymbol{\psi}_i \approx
\sum_{i=1}^\reddim \mathbf{c}^i_j \boldsymbol{\psi}_i\quad
\text{for}\,j = 1, \dotsc, M,
\label{eq:lcomb}
\end{equation}
where $\boldsymbol{\psi}_i$ refers to the $i$-th POD mode. The coefficients
$\mathbf{c}^i_j$ of the linear combination represent the low-dimensional solution
and are usually called \textit{modal coefficients}. Using the matrix
notation, to compute such coefficients it is sufficient a
matrix multiplication $\mathbf{C} = \mathbf{U}_\reddim^T\mathbf{X}$, where the columns
of $\mathbf{C}$ are the vectors $\mathbf{c}^j \in \R^\reddim$ for $j =
1, \dotsc, \reddim$, the matrix $\mathbf{U}_\reddim \in \R^{\fulldim \times \reddim}$ contains the first $\reddim$ POD basis and the superscript $T$ indicates the matrix transpose.

The new pairs $(\mupar_i, \mathbf{c}_i)$, for $i = 1, \dotsc, M$, we
can be exploited in order to find a function $f: \parspace \to \R^\reddim$
capable to predict the modal coefficients for untested parameters. Several
options are available in literature to reach this goal: for instance $n$-dimensional linear
interpolator~\cite{SalmoiraghiScardigliTelibRozza2018, garotta2020}, radial basis functions (RBF)
interpolator~\cite{tezzele2018ecmi}, artificial neural networks~\cite{wang2019non},
Gaussian process regression~\cite{OrtaliDemoRozza2020MINE, guo2018reduced}.
As anticipated, in this work we apply a GPR~\cite{williams2006gaussian}, fitting the distribution of the
modal coefficients with a multivariate Gaussian distribution, such that
\begin{equation}
f(\mupar) \sim \text{GP}(m(\mupar), K(\mupar, \mupar)),
\end{equation}
where $m(\cdot)$ and
$K(\cdot, \cdot)$ indicate the mean and the covariance of the distribution,
respectively. Given a covariance function, an optimization step is required to
set the corresponding hyperparameters. In this contribution we use the squared
exponential covariance defined as $K(x_i, x_j) = \sigma^2
\exp\left(-\frac{\|x_i - x_j\|^2}{2l}\right)$. Once the hyperparameters
($\sigma$ and $l$) of the covariance kernel have been fit to the input dataset,
we can query such distribution to predict the new modal coefficients . Finally
the modal
coefficients are projected back to the high-dimensional vector space
$\R^\fulldim$ using~\eqref{eq:lcomb}.
It is easy to note the differences from the computational point of view between
FOM and ROM: whereas in the full order model it is required to solve a
non-linear problem of dimension $\fulldim$, in the reduced order model to
predict the solution we just need to query a distribution and perform a matrix
multiplication. \RC{From the computational perspective, in fact the cost of the
ROM is mainly due to its construction and not to the prediction phase: relying
on the SVD, the method shows an algorithmic complexity of $\mathcal{O}(\text{min}(\fulldim, M)\,\fulldim M)$.
Thus, dealing with complex FOM as the one presented in this work, POD space
construction can be neglected in the overall computational need.}

On the technical side, we construct and exploit the POD-GPR model using
EZyRB~\cite{DemoTezzeleRozza2018EZyRB}, an open source Python package which deals
with several data-driven model order reduction
techniques, exploiting the library GPy~\cite{gpy} for the GPR implementation.

\section{Optimization procedure with built-in parameters reduction}
\label{sec:asga}

In this work we make use of the active subspaces extension of the
genetic algorithm (ASGA) introduced in~\cite{demo2020asga}. Such
optimization method has been selected as it outperforms standard GA,
especially when high-dimensional target functions are considered.
Its performance have been proved both for classical academic benchmark
functions and for industrial CFD test cases.

The following sections report a description of both the classical
genetic algorithm and the active subspaces technique main features.
Finally, we will discuss how the two algorithms have been combined
to obtain an efficient optimization procedure.

\subsection{Genetic algorithm}
\begin{figure}
\centering
\includegraphics[width=.65\textwidth]{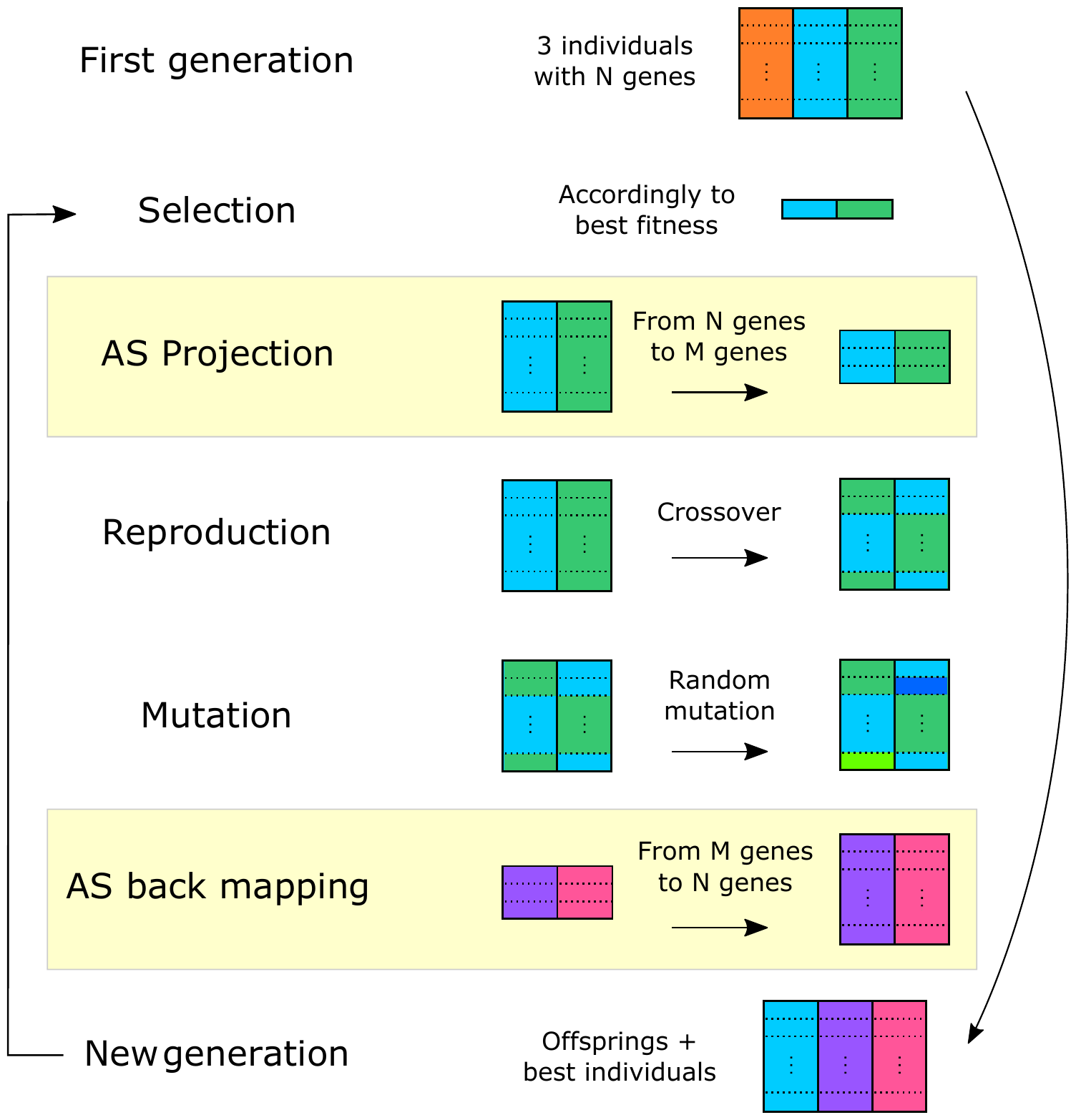}
\caption{Active subspaces-based genetic algorithm scheme. The main
  step of the classical GA are depicted from top to bottom. The yellow
  boxes represent projections onto and from lower dimension active
  subspace. Thus, they are specific to ASGA.\label{fig:asga_scheme}}
\end{figure}

Genetic algorithm (GA) is an optimization algorithm, first introduced
by Holland in~\cite{holland1973genetic}.  Inspired by
natural selection, it falls into the category of population based search
algorithms. For a detailed discussion of the method and its several
modifications we refer the interested reader
to~\cite{katoch2020review, el2006hybrid, sivaraj2011review}. Here,
we briefly present the simplest genetic algorithm, which is composed by three
fundamental steps: \emph{selection}, \emph{reproduction}, and
\emph{mutation}. Such phases are illustrated in Figure~\ref{fig:asga_scheme}
--- which also includes yellow boxes which will be discussed in the
following sections.

The algorithm starts with a random population $\mathcal{S}_0$ composed of $T$
individuals, each one having $r$ genes. In the selection step the
individuals with the best fitness value, for instance $\mathcal{S}_0^{(1)}$
and $\mathcal{S}_0^{(2)}$, are retained. During the reproduction phase,
an offspring $Q$ is produced from these two individuals with a crossover
probability $P_C$. Then, in the last step $Q$ undergoes a mutation with
probability $P_M$, generating $Q^\prime$. This new offspring $Q^\prime$ is
added in the new population $\mathcal{S}_1$ together with the best individuals
of $\mathcal{S}_0$. The three steps are repeated until a predetermined
computation budget is reached.

\subsection{Active subspaces}
The active subspaces (AS)~\cite{constantine2015active, zahm2020gradient,
morhandbook2020} property is an emerging technique for
dimension reduction of parameterized problems. Let us initially assume that
the input/output relationship of the problem under study is represented by
function $f(\mupar) : \Omega \subset \R^n \to \R$.
The reduction is performed by computing a linear
transformation of the original parameters $\mupar_M = \mathbf{A}
\mupar$, in which $\mathbf{A}$ is an $M\times n$ matrix, and $M<n$. 
In the last years AS has been extended to vector-valued output
functions~\cite{zahm2020gradient}, and to nonlinear transformations
of the input parameters using the kernel-based active subspaces (KAS)
method~\cite{romor2020kas}.
AS has been also coupled with reduced order methods such as
POD-Galerkin~\cite{tezzele2018combined} in cardiovascular studies, and
POD with interpolation~\cite{demo2019cras} and dynamic mode
decomposition~\cite{tezzele2020enhancing} for CFD applications.
Application to multi-fidelity approximations of scalar
functions are also presented in~\cite{romor2020pamm, liulin2020}. 

The matrix $\mathbf{A}$ is computed based on the second moment matrix $\mathbf{C}$ of
the target function $f$ gradient. The latter matrix is defined as
\begin{equation}
\mathbf{C} := \mathbb{E}\, [\nabla_{\mupar} f \, \nabla_{\mupar} f
^T] =\int (\nabla_{\mupar} f) ( \nabla_{\mupar} f )^T
\rho \, d \mupar,
\end{equation}
where with $\mathbb{E} [\cdot]$ we denote the expected
value, $\nabla_{\mupar} f \equiv \nabla f({\mupar}) \in \mathbb{R}^n$,
and $\rho : \mathbb{R}^n \to \mathbb{R}^+$ is a probability density
function representing the uncertainty in the input parameters. The
gradients appearing in $C$ are typically approximated~\cite{constantine2015active}
with local linear models, global linear models, GP regression, or finite
difference. The second moment matrix $\mathbf{C}$ is constructed with a
Monte Carlo procedure. We proceed by decomposing the uncentered
covariance matrix as $\mathbf{C} = \mathbf{W} \Lambda \mathbf{W}^T$,
where $\Lambda$ is the diagonal eigenvalues matrix (arranged in
descending order) and $\mathbf{W}$ is the orthogonal matrix containing the
corresponding eigenvectors. To bound the error on the numerical approximation
associated with Monte Carlo simulations, we make use of the gap between the
eigenvalues. Looking at the energy decay, we can select a scalar $M < n$ and
decompose $\Lambda$ and $\mathbf{W}$ as
\begin{equation}
\Lambda =   \begin{bmatrix} \Lambda_1 & \\
& \Lambda_2\end{bmatrix},
\quad
\mathbf{W} = \left [ \mathbf{W}_1 \quad \mathbf{W}_2 \right ],
\quad
\mathbf{W}_1 \in \mathbb{R}^{n\times M},
\end{equation}
where $M$ is the dimension of the active subspace -- which can also be
prescribed a priori. The decomposition described is exploited to map the input 
parameters onto a reduced space. Thus, the principal eigenspace
corresponding to the first $M$ eigenvalue defines the \emph{active subspace}
of dimension $M$. In particular we define the active variable as $\mupar_M := 
\mathbf{W}_1^T\mupar \in \mathbb{R}^M$ and the inactive variable as
$\etapar := \mathbf{W}_2^T\mupar \in \mathbb{R}^{n-M}$.

Exploiting the higher efficiency of most interpolation strategy in lower
dimensional spaces, we can now approximate $f$ using a response surface
over the active subspace, namely
\begin{equation}
g(\mupar_M = \mathbf{W}_1^T\mupar) \approx f(\mupar), \qquad \mupar_M
\in \mathcal{P} := \{ \mathbf{W}_1^T\mupar \; | \; \mupar \in
\Omega \},
\end{equation}
where $\mathcal{P}$ is the polytope in $\mathbb{R}^M$ (the ranges of
the parameters are intervals) defined by the AS. 

The active subspaces technique and several other methods for
parameter spaces reduction are implemented in the ATHENA\footnote{Freely
  available at \url{https://github.com/mathLab/ATHENA}.} Python
package~\cite{romor2020athena}.

\subsection{Active subspaces-based genetic algorithm}
We enhance the classical GA by adding two fundamental steps before the
reproduction and after the mutation phase. These involve the
application of the projection of the current population onto its
active subspace, given a prescribed dimension. So, the idea is to
perform the crossover and the random mutation in the smaller dimension
space. Such space in fact only includes the directions in which the highest
variation of the fitness function $f$ is observed.

By a mathematical standpoint, we add the following
operations to the GA: let $\mathbf{W}_1$ be the eigenvectors defining the
active subspace of the current population, say $\mathcal{S}_0$. We
project its best individuals onto the current active subspace with
\begin{equation}\label{eq:AS_fw_mapping}
s_0^{(1)} = \mathbf{W}_1^T \mathcal{S}_0^{(1)}, \qquad
s_0^{(2)} = \mathbf{W}_1^T \mathcal{S}_0^{(2)},
\end{equation}
where $s_0^{(1)}$ and $s_0^{(2)}$ are the reduced individuals. The
reproduction and mutation steps are performed as usual. The only
difference is that in the described framework they conveniently are
carried out within a smaller dimension space, where reduced number
of genes is exploited for speed up purposes. After these phases are
completed, we obtain the offspring $q$ and $q^\prime$, respectively.
Finally, the back mapping from the active subspace to the full space
is performed by sampling the inactive variable $\etapar$ in order to obtain
\begin{equation}\label{eq:AS_bw_mapping}
Q^\prime = \mathbf{W}_1 q^\prime + \mathbf{W}_2 \etapar, \qquad
\text{with } -\mathbf{1} \leq Q^\prime \leq \mathbf{1},
\end{equation}
where $\mathbf{1}$ denotes a vector with all components equal to
$1$ --- the original parameters are usually rescaled in $[-1, 1]^n$
before applying AS ---. We remark that there is in principle the possibility that
multiple points in the full space are mapped onto the same reduced point
in the active subspace. Hence, the number $B$ of individuals resulting
from the back mapping is an hyperparameter which can be prescribed a
priori. For the specifics about this procedure please
refer to~\cite{demo2020asga}. In Figure~\ref{fig:asga_scheme} we emphasized
with yellow boxes the new fundamental steps represented by
Equations~\eqref{eq:AS_fw_mapping} and \eqref{eq:AS_fw_mapping}. For
the actual implementation of the genetic algorithm part we used DEAP~\cite{fortin2012deap}.

\section{Numerical results}
\label{sec:results}

\begin{figure}[t]
\centering
\includegraphics[width=1.\textwidth]{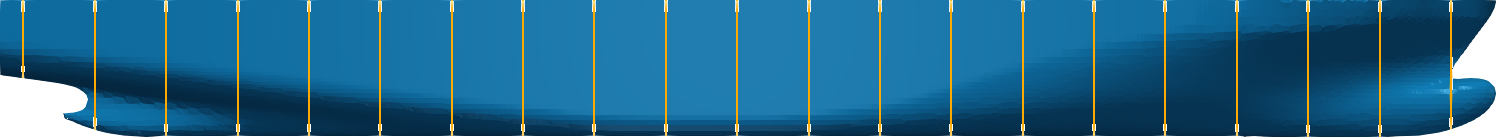}
	\caption{The surface of the DTC hull. The highlighted sections
          divide the ship into $20$ equispaced chunks at the
          free-surface level.}\label{fig:dtc}
\end{figure}
In this section, we describe the application of the proposed optimization pipeline
to the DTC hull surface. \RA{Table~\ref{tab:dtc} shows the main particulars in
the design loading condition at model scale (which is set to $1:59.407$)}.
This will provide a test case which closely simulates
a typical workflow for industrial hull design problems. Figure~\ref{fig:dtc}
shows the original CAD geometry of the hull used in this work, where we marked
$21$ longitudinal sections which divide the ship into $20$ equispaced chunks.
Such $21$ slices will be referred to as \textit{sections} during the results
discussion, and are numbered from $1$ to $21$ going from the ship stern to its bow.
\begin{table}[b]
\centering
    \caption{\RA{Main quantities of the DTC at scale model.}}\label{tab:dtc}
    \begin{tabular}{llc}
    \hline
        \textbf{Quantity} & \textbf{Value}\\
    \hline
        Length between perpendiculars $L_{pp}~[\si{\meter}]$ & $5.976$ \\
        Waterline breadth $B_{wl}~[\si{\meter}]$ & $0.859$ \\
        Draught midships $T_{m}~[\si{\meter}]$ & $0.244$ \\
        Volume displacement $V~[\si{\meter^3}]$ & $0.827$ \\
        Block coefficient $C_B$ & $0.661$\\
    \hline
    \end{tabular}
\end{table}

The structure of this section mirrors that of the whole article,
reporting the intermediate results of all the methods employed throughout the
optimization pipeline.

\subsection{Self-learning mesh morphing parameters}

To set up the FFD hull surface deformation, we position the control points
lattice in order to control the immersed part of the ship prow region. The
equispaced control points are positioned as follows:
\begin{itemize}
\item $\boldsymbol{x}$ {\bf axis:} 7 points layers located on sections 10, 12, 14, 16, 18, 20, 22;
\item $\boldsymbol{y}$ {\bf axis:} 11 points layers that cover the whole hull beam,
      with the second and the second-to-last positioned on the lateral walls of the ship;
\item $\boldsymbol{z}$ {\bf axis:} 7 points layers that cover the whole hull draft, aligning
      the $2^{nd}$ and the $5^{th}$ of them to the hull bottom and to the
      waterline, respectively.
\end{itemize}
As can be appreciated by the values reported, to distribute the FFD control points,
we have made use of an additional 22nd \emph{virtual} section located ahead of the
bow. The motion of the $7\times 11\times 7 = 539$ points is governed by only 10
parameters, which are described in Table~\ref{tab:ffd}. We point out that the displacement
of all the boundary points in the $x$ and $z$ direction is set to zero so as to enforce
surface continuity. In addition, the displacement of the points on the internal $x$ and $z$
layers closest to the boundary ones is also set to zero so as to enforce continuity of
all surface derivatives. Finally, the hull symmetry along $y$ direction is ensured by
selecting symmetric values for parameters associated to $x$ and $z$ displacements, as
well as antisymmetric values for parameters associated to $y$ displacements (the latter
points are also indicated in the table by the corresponding footnote).

\begin{figure}[b]
\centering
\includegraphics[width=.48\textwidth]{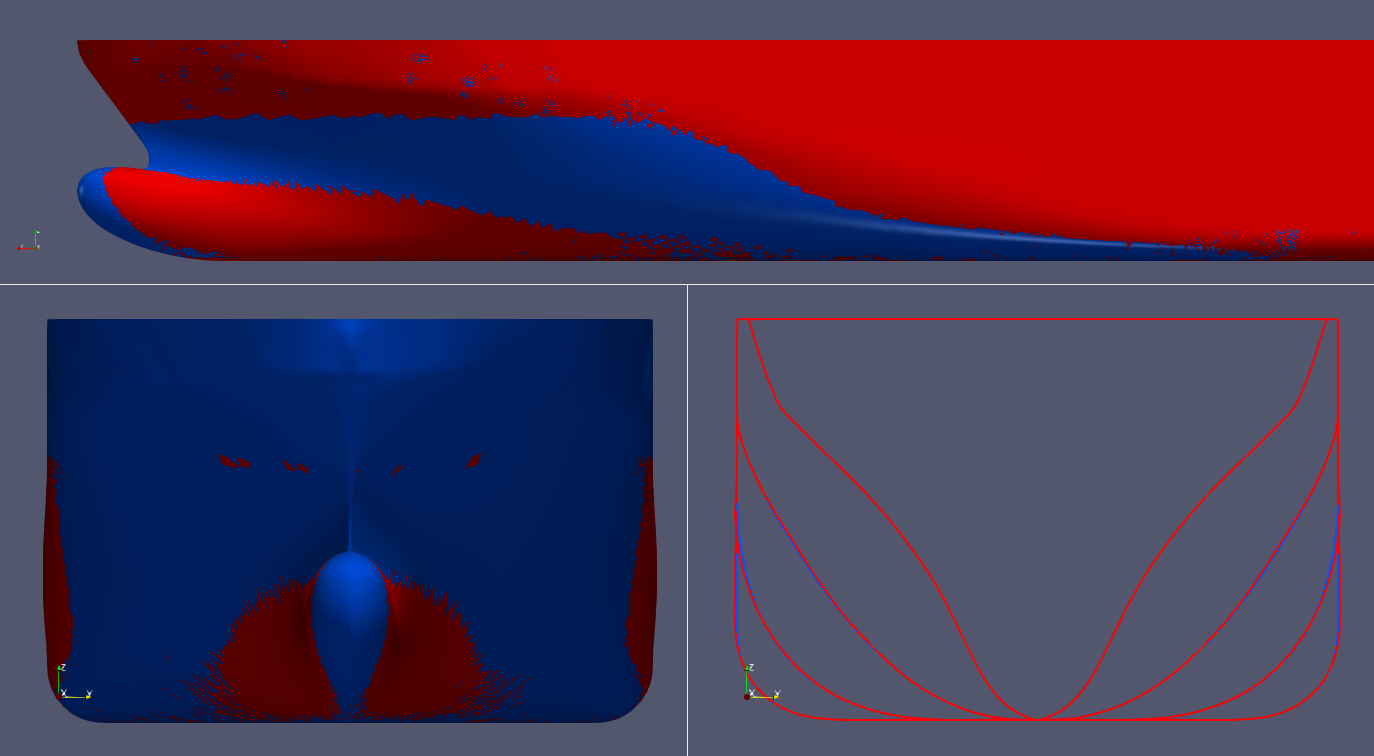}
\includegraphics[width=.48\textwidth]{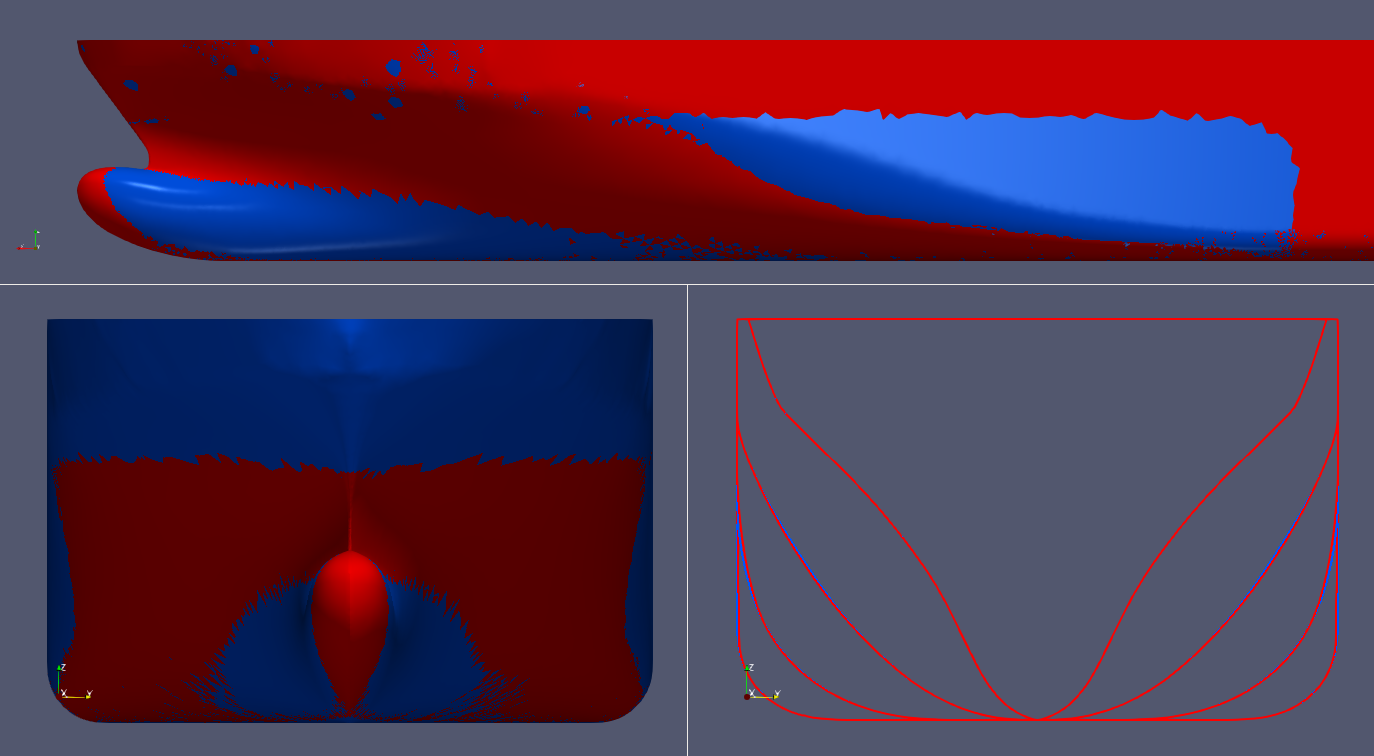}
	\caption{Visual examples of hull deformation with $\mupar =
          [-0.2]^{10}$ (on left) and $\mupar = [0.2]^{10}$ (on
          right). The red surface refers to the deformed ships, while
          the blue one is the original hull.\label{fig:ffd_deform}} 
\end{figure}

Once defined the geometric parameters $\mupar = [\mu_0, \dotsc ,
\mu_9 ]$, we set the parametric space to $\parspace = [-0.2,
0.2]^{10}$. The parameter space boundary values are selected so as to
obtain feasible deformations from an engineering point of
view and, at same time, to explore a large variety of possible
shapes. Figure~\ref{fig:ffd_deform} shows the two ``extreme'' hull
deformations, obtained setting all the parameters equal to the lower
and upper bound of the space, respectively.

\begin{table}[t]
\centering
\caption{FFD control points displacement. The indices refer to the relative
	position of the points within the lattice.The layers order, which starts
	from 0, is maintained consistent with the reference system. The intervals
	indicated by the -- symbol are inclusive.}\label{tab:ffd}
\begin{tabular}{ccccc}
\hline\hline
    \multicolumn{3}{c}{Lattice Points} & \multirow{2}{*}{Parameter} & \multirow{2}{*}{Displacement direction}\\
    index $x$ & index $y$ & index $z$ & & \\
\hline\hline
\rowcolor{lightgray}
    2 & 0 & 2--4 & $\mu_0$ & $x$\\
    2 & 10 & 2--4 & $\mu_0$ & $x$\\
\rowcolor{lightgray}
	3 & 0 & 2--4 & $\mu_1$ & $x$\\
	3 & 10 & 2--4 & $\mu_1$ & $x$\\
\rowcolor{lightgray}
	4 & 0 & 2--4 & $\mu_2$ & $x$\\
	4 & 10 & 2--4 & $\mu_2$ & $x$\\
\rowcolor{lightgray}
	4 & 2--4  & 2 & $\mu_3$ & $y$\\
	4 & 6--8 & 2 & $-\mu_3$ \footnotemark[2] & $y$\\
\rowcolor{lightgray}
	4 & 2--4 & 3 & $\mu_4$ & $y$\\
	4 & 6--8 & 3 & $-\mu_4$ \footnotemark[2] & $y$\\
\rowcolor{lightgray}
	4 & 2--4 & 4 & $\mu_5$ & $y$\\
	4 & 6--8 & 4 & $-\mu_5$ \footnotemark[2] & $y$\\
\rowcolor{lightgray}
	3 & 2--4 & 2 & $\mu_6$ & $y$\\
	3 & 6--8 & 2 & $-\mu_6$ \footnotemark[2] & $y$\\
\rowcolor{lightgray}
	5 & 2--4 & 3 & $\mu_7$ & $y$\\
	5 & 6--8 & 3 & $-\mu_7$ \footnotemark[2] & $y$\\
\rowcolor{lightgray}
	4 & 0--1 & 2 & $\mu_8$ & $z$\\
	4 & 9--10 & 2 & $\mu_8$ & $z$\\
\rowcolor{lightgray}
	5 & 0    & 3 & $\mu_9$ & $z$\\
	5 & 10    & 3 & $\mu_9$ & $z$\\
\hline\hline
\end{tabular}
\end{table}


The FFD deformation of the hull points has been extended to the
nodes of the volumetric grid for the CFD simulations making use of
the Beckert-Wendland radial basis function
kernel~\cite{beckert2001multivariate}, defined as follows
\begin{equation}
\label{eq:Beck-Wendl}
\varphi_j(||\boldsymbol{x}-\boldsymbol{x}_j||)=\left(1-\frac{||\boldsymbol{x}-\boldsymbol{x}_j||}{R}\right)^4_+
                                               \left(1+4\frac{||\boldsymbol{x}-\boldsymbol{x}_j||}{R}\right),
\end{equation}
where $R>0$ is a prescribed finite radius and the $(\cdot)_+$ symbol indicates the positive part. 

\begin{figure}[h]
\centering
\includegraphics[width=1.\textwidth]{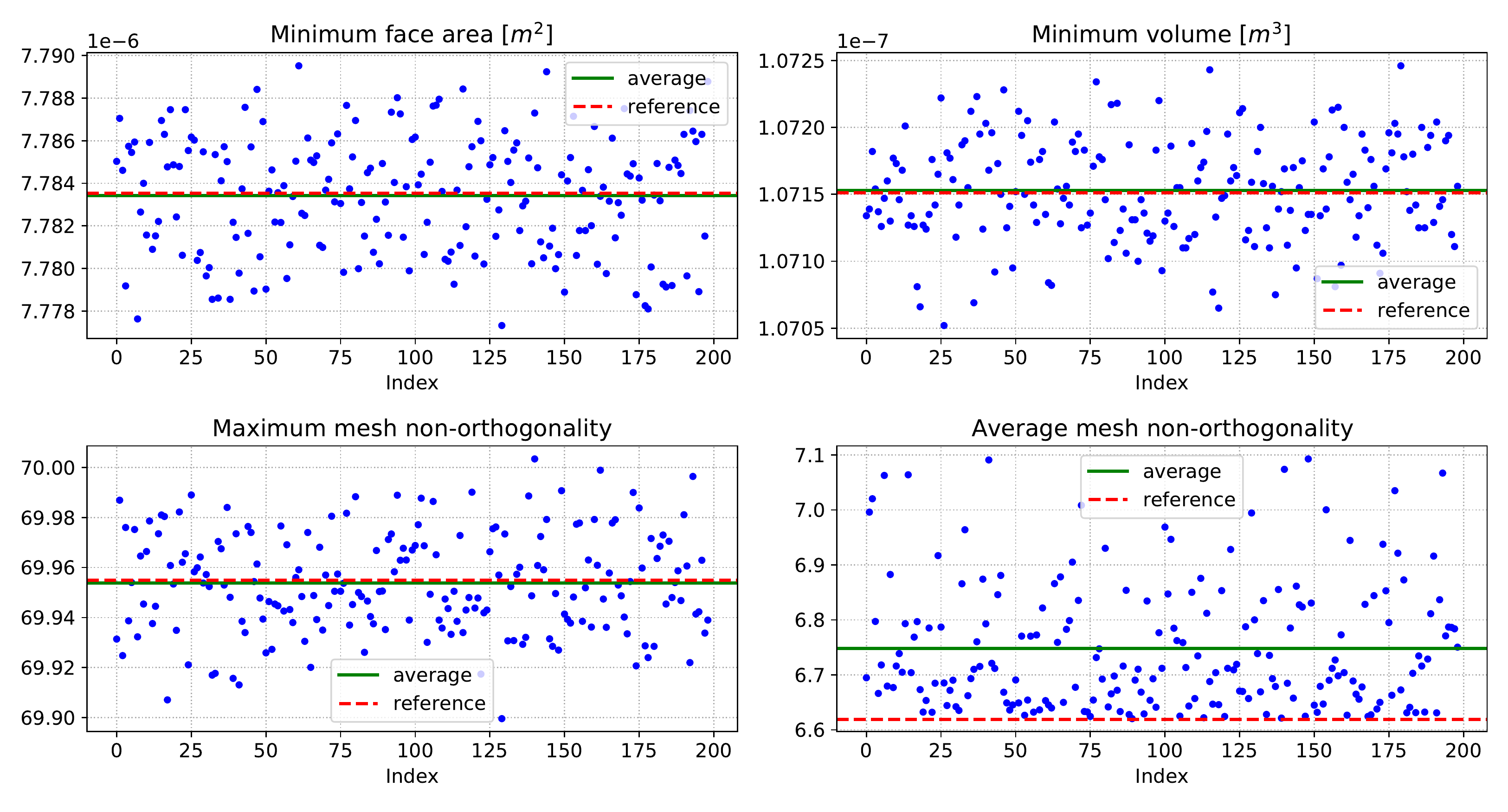}
\caption{Values of the main mesh quality indicators as reported by checkMesh
         utility of OpenFOAM library\RB{, as a function of the index corresponding to each of the $200$
         volumetric meshes produced for the simulation campaign}. \label{fig:mesh_quality}}
\end{figure}

The output of the OpenFOAM library checkMesh utility has been used to
assess the quality of the grids obtained with the combined FFD/RBF methodology.
Figure~\ref{fig:mesh_quality} presents some of the main quality indicators of
the $200$ meshes generated for the present campaign, as computed by checkMesh. In
particular, the indicators considered are minimum face area (top left plot),
minimum cell volume (top right plot), maximum mesh non-orthogonality (bottom left plot)
and average mesh non-orthogonality (bottom right plot). In all the diagrams, the
vertical axis refers to the mesh quality indicator considered, while the variable
associated with the horizontal axis is the index corresponding to each of the$ 200$
volumetric meshes produced for the simulation campaign.

\footnotetext[2]{Imposed for $y$ symmetry conservation.}

The minimum face area and minimum cell volume results indicate that the morphing
procedure does not produce negative cells or faces which would impair the simulations. In fact,
the average of both indicators across the $200$ grids produced is extremely close to the
corresponding value of the original grid. The lowest value of minimum face area observed 
in the $200$ grids generated is less than $0.1\%$ off the original value, while the lowest value
of minimum cell volume observed is merely $0.01\%$ off the original mesh minimum cell volume.
Such trend is confirmed by the maximum non-orthogonality values reported in the bottom
left diagram. In the plot, is possible to appreciate that the average over the $200$ grids
produced falls exactly on value of the original mesh, and the highest difference with
respect to the original mesh non-orthogonality is merely $0.05\%$.
These values ensured that all the simulations in the present campaign could be completed
in fully automated fashion without crashes were reported or significant issues were observed.
The results reported in the bottom right plot indicate that the effect of the mesh
morphing algorithm proposed is that of increasing the grid average non-orthogonality
values. This is somewhat expected, as the original volumetric grid in this work
was generated making use of the snappyHexMesh tool of the OpenFOAM library. In such 
framework, most of the cells in the internal regions of the domain are substantially
the result of an octree refinement of an original block mesh aligned with the
coordinate axes. It is clear that the RBF procedure described in
Section~\ref{sec:shape} does quite clearly alter in a non negligible way the orthogonal
angles of a portion of the hexahedral cells produced by snappyHexMesh. Yet, the average 
increase in the average mesh non-orthogonality index is $2\%$, while the maximum increase
observed is $7.2\%$, which are values that should not significantly affect the results
of the simulations.

\subsection{Reduced order model construction}
We set the full order model in scale $1:59.407$, keeping it unaltered from the
original work mainly for validation purpose. The computational
domain\RB{, that is a parallelepiped of dimension $[ -26, 16 ] \times [ -19, 0 ] \times [ -16, 4 ]$ along $x, y$ and $z$ directions} is
discretized in $8.5 \times 10^5$ cells, with anisotropic vertical refinements
located particular in the free-surface region, in order to avoid a too diffusive
treatment of the VOF variable. 
\RB{Boundaries of such domain are imposed as follows:
\begin{itemize}
\item at the \emph{inlet} we set constant velocity, fixed flux condition for
    the pressure and a fixed profile for the VOF variable;
\item at the \emph{outlet} we set constant average velocity, zero-gradient
    condition for the pressure and variable height flow rate condition for VOF
    variable;
\item at the bottom and lateral planes, we impose symmetric conditions for all
    the quantities;
\item at the top plane, we set a pressure inlet outlet velocity condition for
    the velocity and nil pressure; VOF variable is fixed to 1 (air);
\item at the hull surface, we impose no-slip condition for velocity, fixed flux
    condition for the pressure and zero-gradient condition for VOF variable.
\end{itemize}
The adopted solver is \emph{interFoam}, which is able to solve the Navier
Stokes equations for two incompressible, isothermal immiscible fluids. Time
discretization uses a first order implicit scheme with local-step, since we are
interested to the steady solution. For the spatial discretization, we apply a
Gaussian integration using second order upwind scheme for divergence operators
and linear interpolation for gradient and laplacian operator. 
By imposing a inlet velocity of $1.668~\si{\meter/\second}$, the Froude number
is around $0.22$.
}
The time required to converge to the steady
solution within such setting on a parallel machine ($32$ processors) is
approximately $2$ hours.

For the construction of the reduced order model, we randomly sample the
parametric space with uniform distribution. We performed 203 simulations with
the full order model, collecting the corresponding pressure and shear stress
distributions (the latter implicitly containing the distribution of the VOF variable)
over the hull surface. Thus, only the surface fields are
considered at the reduced level.
We then flatten the shear stress vector field in order to construct two
snapshots matrices, one for the pressure and one for the stress. Both are then
decomposed using POD technique. The number of modes considered is fixed to $20$.
Approximating the manifold with the GPR method, we obtain two different POD-GPR
model that approximate the pressure field and the shear stress field. Such
quantities are used for the computation of the objective function
during the optimization procedure.

Even if the difference of hardware used for full order model simulations and
for reduced order approximation limits the possible speedup obtained --- 
a HPC facilities versus an ordinary personal computer ---, we achieve
satisfactory computational gain. In fact, whereas the FOM lasts approximately two
hours, the ROM approximation only consisting in two distribution queries and two matrix
multiplications, takes less than $1$ second in a single-processor environment. 
Such results are very effective in the framework of an iterative process, as
the optimization pipeline here proposed. The overall time is in fact mainly
constituted by the initial FOM simulations needed for the offline database,
while the ROM approximation can be considered negligible from the computational
point of view. Moreover, it can be performed on significantly less powerful machines.

Adopting data-driven methodologies rather than projection-based ones has
different advantages which we have already discussed, but shows also some
drawback in the error bounding.
For an a posteriori quantification of the ROM accuracy we need then to validate
the approximated optimal result by carrying out a FOM simulation. We remark that
we consider the output of such simulation as truth solution. This requires an
additional computational cost, but allow also for an effective refinement of
the ROM. Once a geometrical configuration is validated in such fashion, depending
on the error observed we can add this last snapshot to the database and re-build
the ROMs.

\subsection{Optimization procedure}
We first define the objective function we applied to the optimization procedure.
The quantity to minimize is \RB{the total resistance coefficient $C_t$}, which is defined as
\begin{equation}
\label{eq:obj}
\underset{\mupar}{\text{min}}\;C_t\equiv
\underset{\mupar}{\text{min}}
	\int_{\Omega(\mupar)} \frac{\tau_x \rho - p n_x}{\frac{1}{2}\rho V^2 S},
\end{equation}
where $\tau_x$ is the $x$-component of the shear stress, $\rho$ is the fluid
density, $p$ indicates the pressure, $n_x$ the $x$-component of the surface
normal, $V$ and $S=\Delta^{2/3}$ the reference fluid velocity and the reference surface,
respectively. As reported, the CFD simulations have been carried out in fixed sink and
trim conditions. Thus, the specific reference surface used to obtain $C_t$ has been
selected to penalize hulls obtaining resistance gains through immersed volume
reduction. All the geometrical quantities, as well as the normals and the reference
surface depend by the imposed deformation. Thus, to evaluate the $C_t$ for
any design, we deform the hull surface using the FFD map, then
project the ROM approximated fields --- pressure and shear stress --- on it to
numerically compute the integral defined in Equation~\eqref{eq:obj}.

Regarding the ASGA hyperparameters, 
we set the probability of crossover and mutation as $P_C = P_M =
0.5$. For each solutions database we perform an optimization run with
ASGA composed by $150$ generations, with an initial random population
of $100$ individuals and an offspring of $20$ individuals. The number
of points returned by the AS back mapping is $B=2$, while the
dimension of the AS is set to $1$ for every population. The covariance
matrix for the active subspace computation is approximated using local
linear models~\cite{constantine2015active}.
\begin{figure}[h]
\centering
\includegraphics[width=.85\textwidth]{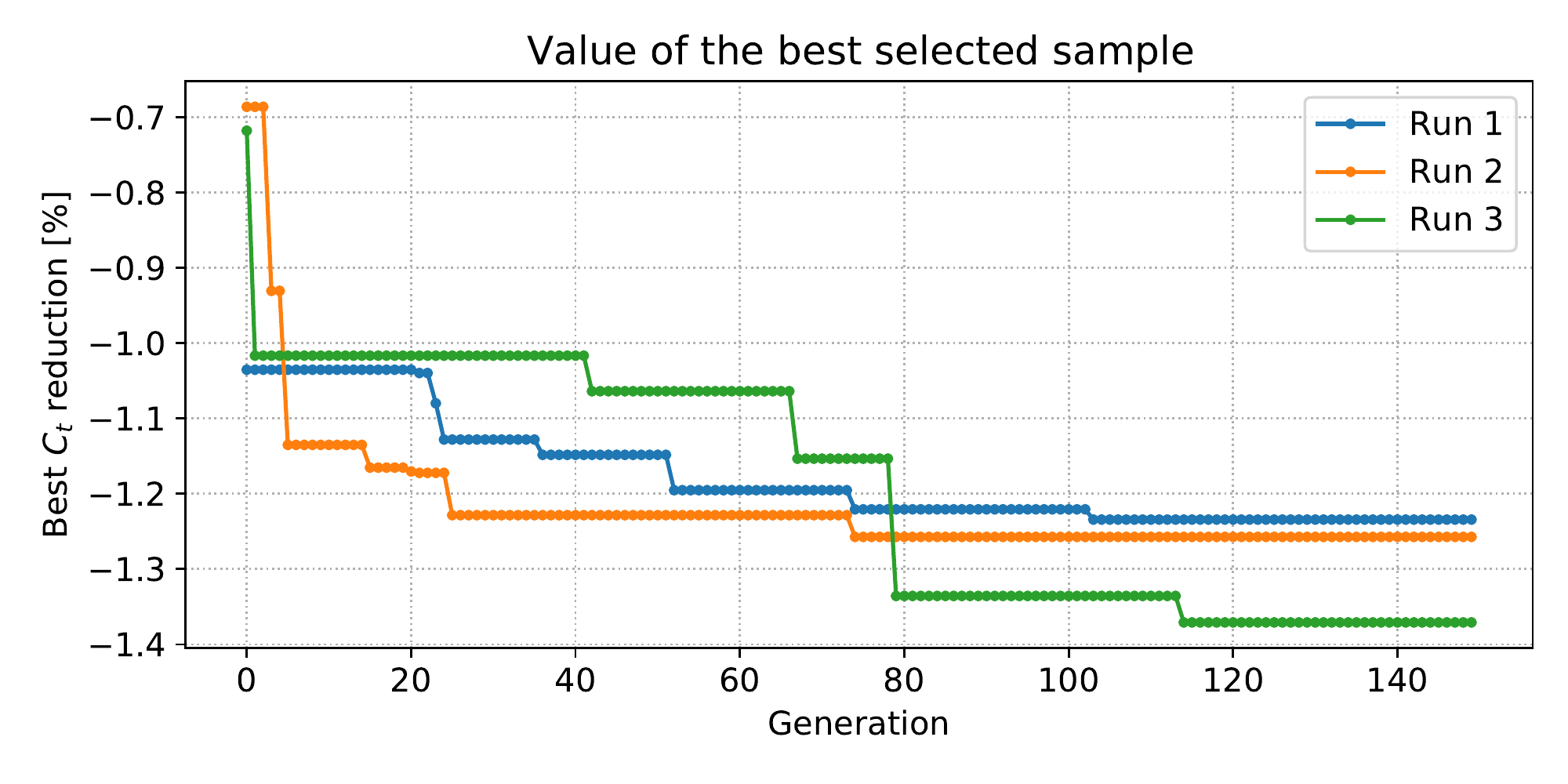}
\caption{ASGA runs. The reduction of the $C_t$ is to be intended with
  respect to the undeformed reference hull.\label{fig:asga_runs}} 
\end{figure}

For each optimum
found by ASGA we run a new high-fidelity simulation for validating the
approximated $C_t$, adding the high-fidelity snapshots to the database in order
to refine the POD-GPR model. In Figure~\ref{fig:asga_runs} we show the
comparison
of all the runs. The third and last optimization reached a reduction of $\sim
1.4\%$ of the $C_t$ coefficient compared to the original shape.

\begin{figure}[bp]
\centering
\includegraphics[width=.65\textwidth]{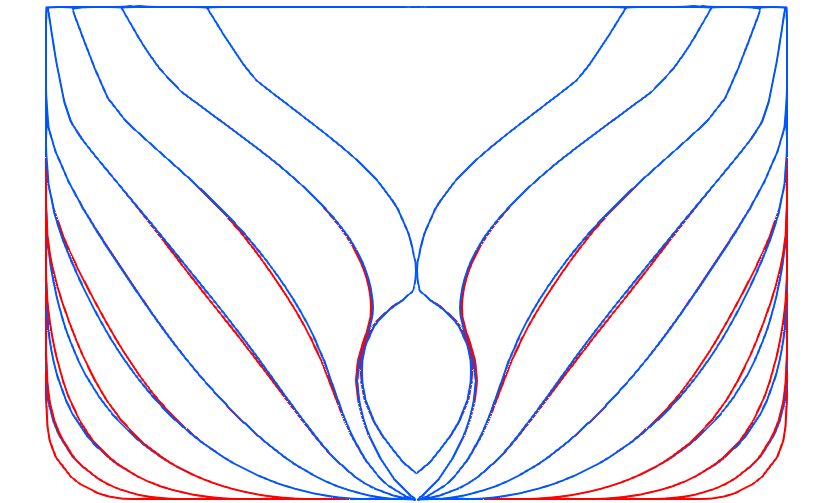}
	\caption{The sections (from $10$ to $20$) of the original ship in blue and of the optimized one in~red.}\label{fig:opt}
\end{figure}

\begin{figure}
\centering
\includegraphics[width=.9\textwidth]{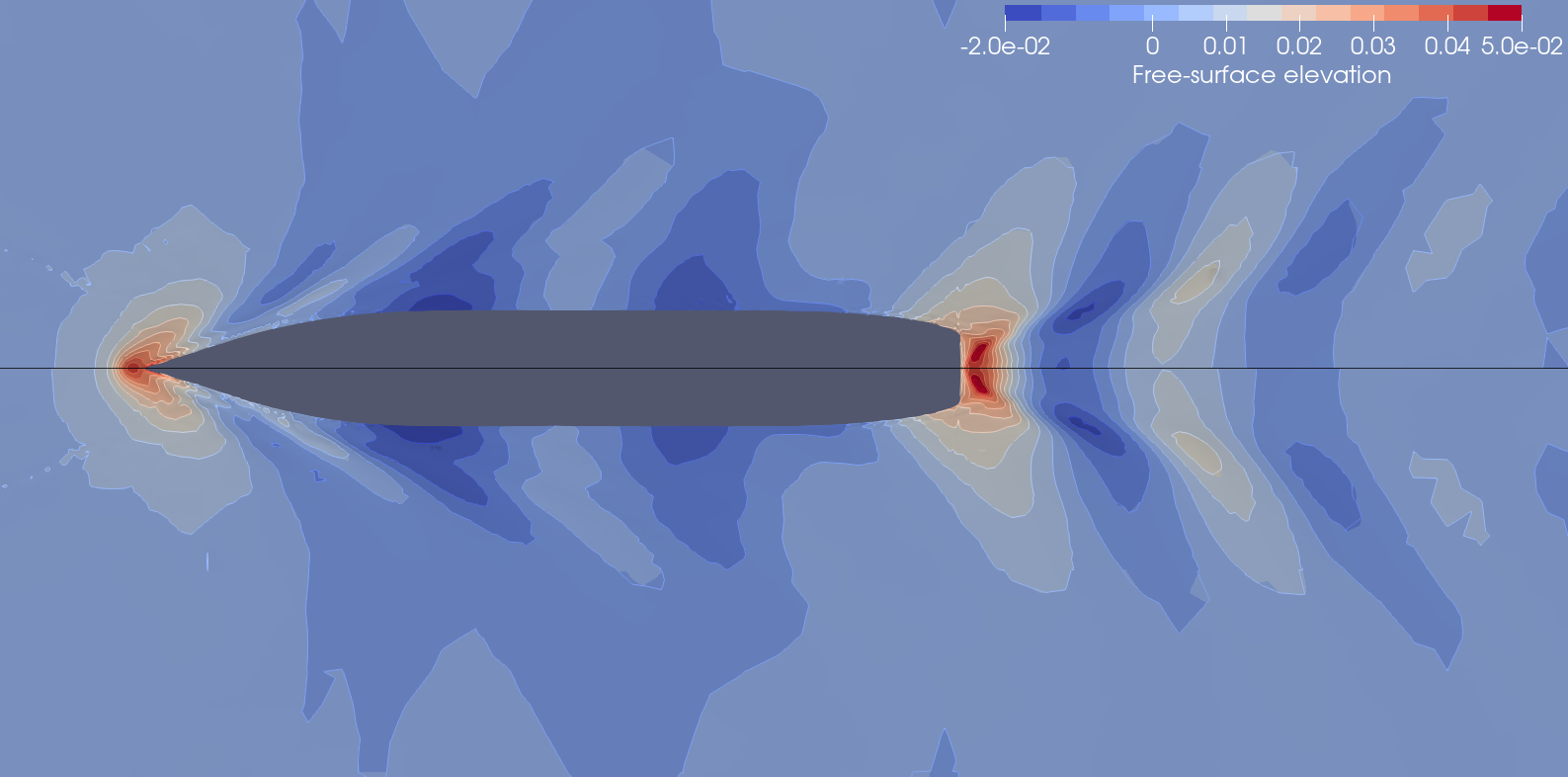}
\caption{Contours of free surface elevation field around the
         original hull (top half)
  and optimal (bottom half).}\label{fig:elevations} 
\end{figure}

Figure~\ref{fig:opt} presents the
frontal sections of the optimal shape compared to the undeformed one, showing a
volumetric increment in the frontal part which balances the reduction near the
central zone. The a posteriori validation
confirmed the positive trend: the $C_t$ coefficient of the optimal shape is
$1.2\%$ less, with a relative error of the ROM model of $0.18\%$. 
\RB{
As is appreciable in Figure~\ref{fig:opt}, the
optimal hull has a wider section in the region immediately downstream with respect
to the bulbous bow, while it appears slightly narrower in the middle ship sections.
The immersed volume of the optimized hull is only $0.08\%$ different from that of the
original hull, which suggests that the $C_t$ reduction obtained is the result
of a total resistance reduction. A possible interpretation of such a resistance decrease
is that having a more streamlined hull along the longitudinal direction, is likely
able to reduce the extent and dimension of the separation bubble located on the
side of the bulbous bow, and corresponding to the dark blue strip visible
in the wall shear stress contours presented in Figures~\ref{fig:org_wss} and
\ref{fig:opt_wss}. As a consequence, the optimal hull presents slightly lower
pressures with respect to the original hull, in the region located downstream of
the bulbous bow. Such a minimal reduction is hardly noticeable in the
pressure contour plots presented in Figures~\ref{fig:org_press} and~\ref{fig:opt_press}.
More appreciable differences are visible instead in the free surface elevation plot presented in
Figure~\ref{fig:elevations}. Reducing the extent of the aforementioned detachment bubble, the
shape modification leading to the optimal hull has the effect of moving forward the
trough which follows the bow. This indicates that the pressures in the bow region are
reduced, which results in a net decrease of the resistance pressure component. In fact,
this leads to a $4.92\%$ reduction in the pressure component of the
resistance, against a more modest $0.55\%$ reduction of viscous resistance. Yet,
considering that the latter component accounts for approximately $83\%$ of the total
resistance, this translates into the $1.2\%$ reduction reported. Finally, to exclude the
possibility that the differences observed in the total resistance coefficient values
are a result of possible discretization error due to the mesh morphing procedure,
we report that the average and maximum values of wall $y^+$ of the optimized hull
do not significantly differ from those obtained with the original one. The
average and maximum wall $y+$ values for the original hull simulation are $6.18426$ and $99.5631$,
respectively, while the corresponding average and maximum values for the optimized hull
are $6.19071$ and $99.6255$, respectively.} We point out that the $y^+$ maxima here
reported for the DTC tutorial appear outside of the range prescribed for the turbulence
model here used. Yet, the accuracy of the DTC tutorial results suggests that
maxima $y^+$ is likely located outside the water. In fact, considering the
small density of air with respect to water, the impact of the resulting inaccurate
estimation of surface derivatives is minimal.

\RA{
We remark that the POD-GPR model approximates the distribution of the output of
interest, not the objective function --- which is computed using the predicted
fields. For this reason, we can also compare the pressure and shear stresses
over the optimal hull with respect to the undeformed one.
Figures~\ref{fig:org_wss} and~\ref{fig:org_press}
 present the graphical investigations about the ROM approximation
error distribution over the undeformed hull, both for pressure and stresses
distributions.  For a more
realistic comparison, we specify that the FOM snapshots referring to the
undeformed geometry has been removed from the database, emulating the
approximation any untested parameter. We proceed in the same way also for the
optimal shape (Figures~\ref{fig:opt_wss} and~\ref{fig:opt_press}), not only to
measure the accuracy of the POD-GPR model, but also for investigating the
reasons of the $C_t$ reduction from a physical perspective. The absolute error
is quite small, but it is possible to note that for both the fields it is
mainly concentrated along the free-surface.

Comparing the original hull with the optimal one we emphasize that the
optimal shape seems to be able to slightly reduce the height of the
wave created by its body, inducing a reduction of the 
wet surface. The friction resistance computed as the integral of the $x$
component of shear stresses over the two hulls shows in fact this marginal gain:
the $12.76~\si{\newton}$ of the original ship becomes $12.69~\si{\newton}$ in
the optimal configuration. However, the main contribution of the resistance
reduction comes from the pressure resistance. While in the original shape we
measure $2.64~\si{\newton}$, in the optimized such quantity decreases to
$2.51~\si{\newton}$.

\begin{figure}[h]
\centering
\includegraphics[width=1.\textwidth]{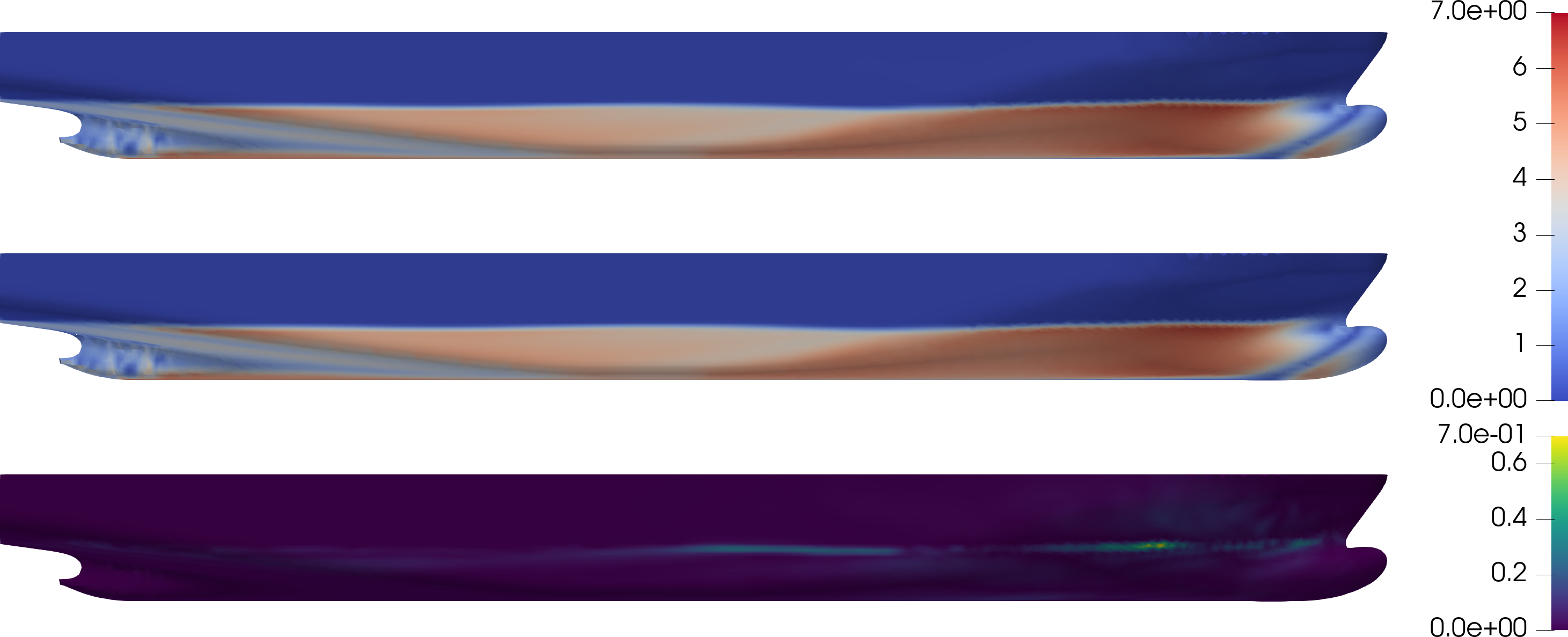}
    \caption{\RA{Distribution of the shear stresses measured in Pascal
      over the undeformed hull: the FOM validation (top) is compared
      to the ROM approximation (middle) and the absolute error is
      shown (bottom).}}\label{fig:org_wss} 
\end{figure}
\begin{figure}
\centering
\includegraphics[width=1.\textwidth]{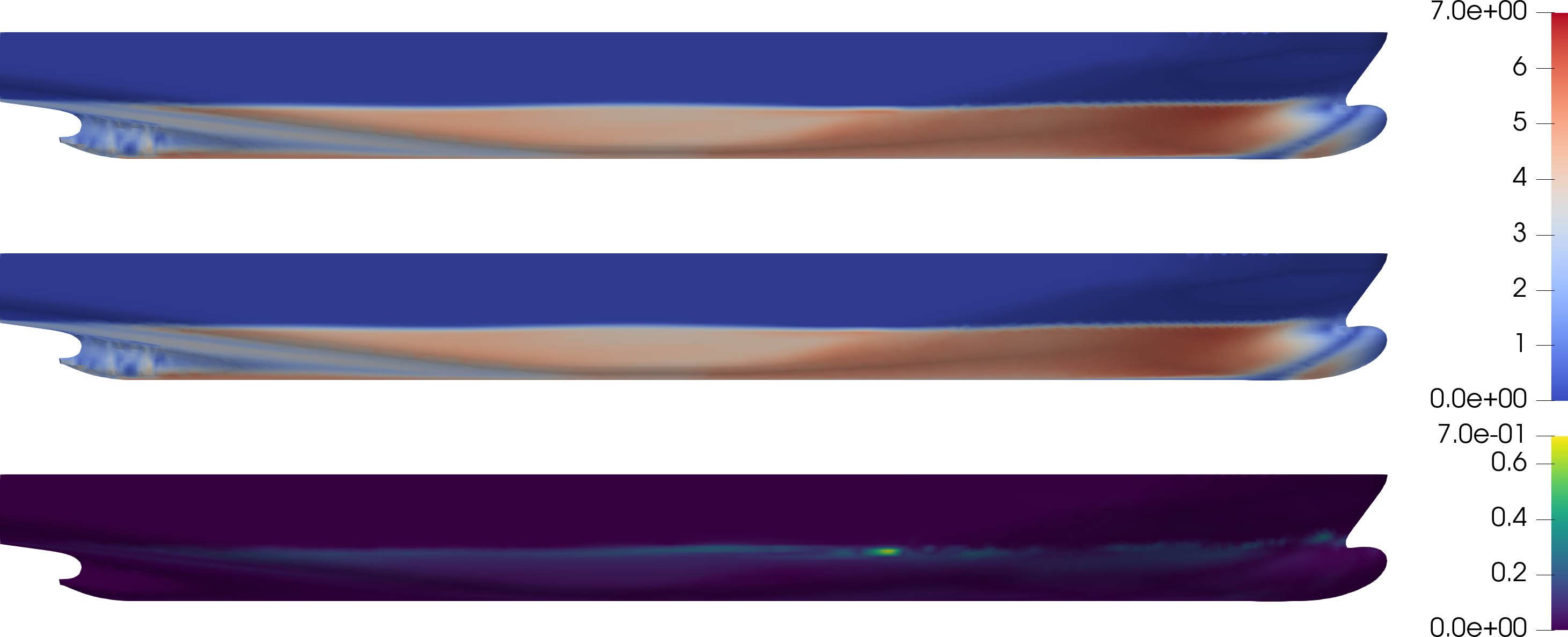}
\caption{\RA{Distribution of the shear stresses measured in Pascal over
  the optimal hull: the FOM validation (top) is compared to the ROM
  approximation (middle) and the absolute error is shown
  (bottom).}}\label{fig:opt_wss} 
\end{figure}
\begin{figure}
\centering
\includegraphics[width=1.\textwidth]{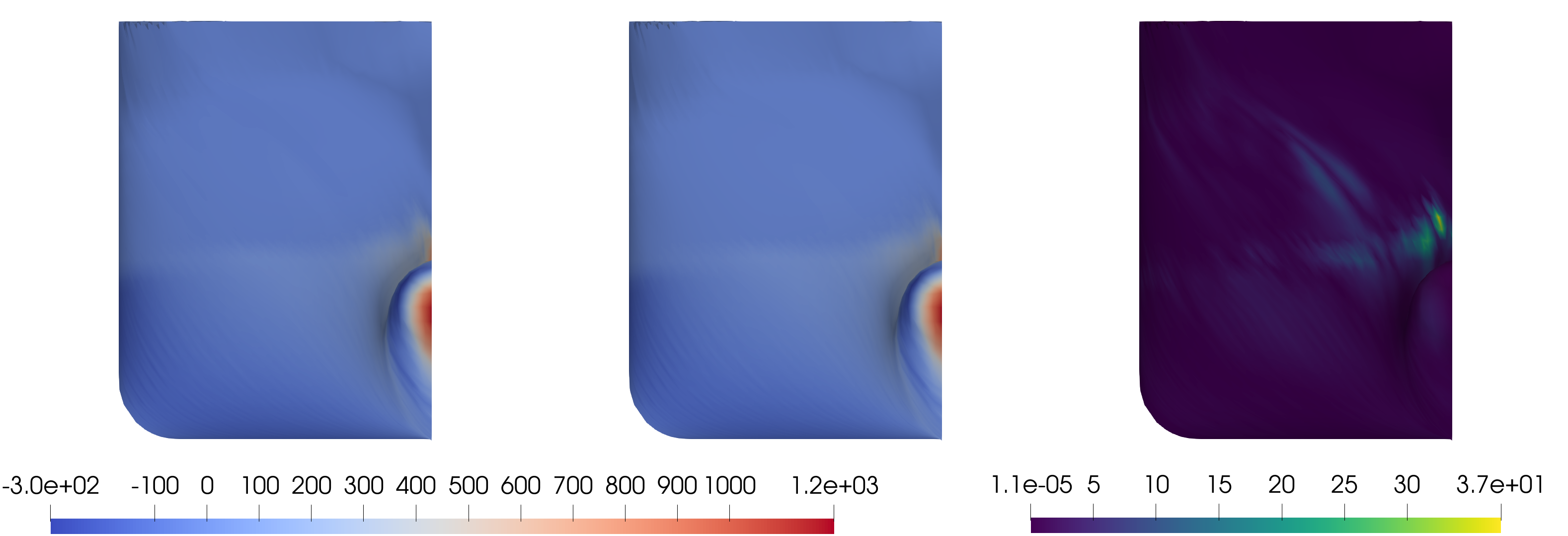}
\caption{\RA{Distribution of pressure measured in Pascal over the
  undeformed hull: the FOM validation (left) is compared to the ROM
  approximation (center) and the absolute error is shown
  (right).}}\label{fig:org_press} 
\end{figure}
\begin{figure}
\centering
\includegraphics[width=1.\textwidth]{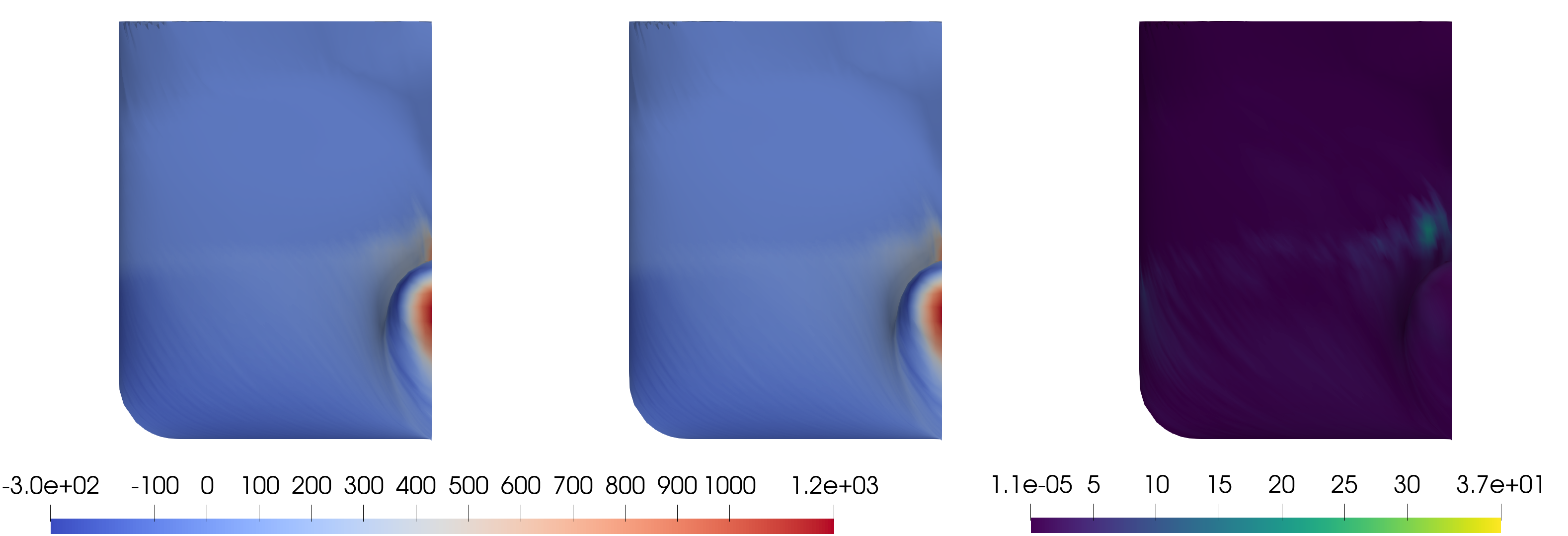}
\caption{\RA{Distribution of the pressure measured in Pascal over the
  optimal hull: the FOM validation (left) is compared to the ROM
  approximation (center) and the absolute error is shown
  (right).}}\label{fig:opt_press} 
\end{figure}
}

\section{Conclusions}
\label{sec:conclusions}

In this work we presented a complete numerical pipeline for the hull
shape design optimization of the DTC benchmark hull. We proposed a
self-learning geometrical deformation 
technique, where different morphing methods are coupled together to
propagate surface deformations to volumetric meshes. Though in this
work we used a FFD approach for the CAD modifications, we emphasize
that our methodology can exploit any surface deformation.
The optimization procedure is based on a coupling between active
subspaces and genetic algorithm, called ASGA. For the evaluation of
the \RB{total} resistance coefficient for new untested parameters we exploits the
non-intrusive data driven reduced order method called POD-GPR. This
results in a great computational saving for the computation of the
pressure and viscous forces fields, while preserving a good accuracy.
We performed $3$ optimization runs, with high-fidelity
validation of the approximated optimum and enrichment of the solutions
database to increase the accuracy of the ROM in its neighborhood.
We obtained a reduction of the \RB{total} resistance coefficient equal to $1.2\%$
with respect to the original reference hull.

In the future, further investigations will be carried out to study a
dynamic selection of the active subspace dimension, and a varying
number of points returned by the back mapping procedure. Further
improvements in the shape parameterization algorithms could be obtained
improving the efficiency of the RBF weights computation. This could
be obtained with a smarter selection of the RBF control points or,
in a more invasive fashion, by resorting to fast algorithms --- such
as Fast Multipole Method --- for the computation of the
control points mutual distances.

\section*{Acknowledgements}
This work was partially supported by an industrial Ph.D. grant
sponsored by Fincantieri S.p.A., and partially funded by the project UBE2 -
``Underwater blue efficiency 2'' funded by Regione FVG, POR-FESR 2014-2020,
Piano Operativo Regionale Fondo Europeo per lo Sviluppo Regionale.  It was also
partially supported by European Union Funding for Research and Innovation ---
Horizon 2020 Program --- in the framework of European Research Council
Executive Agency: H2020 ERC CoG 2015 AROMA-CFD project 681447 ``Advanced
Reduced Order Methods with Applications in Computational Fluid Dynamics'' P.I.
Gianluigi Rozza.

\bibliographystyle{abbrv}

\end{document}